\documentclass[11pt]{amsart}
\usepackage{amssymb}
\usepackage[matrix,arrow,curve]{xy}
\usepackage{graphics}
\usepackage{euscript}

\newcommand{\R}{\mathbb{R}}
\newcommand{\Z}{\mathbb{Z}}

\newcommand{\C}{\mathbb{C}}

\newcommand{\iso}{\cong}           
\newcommand{\htp}{\simeq}          
\newcommand{\CP}[1]{\C {\mathrm P}^{#1}}

\theoremstyle{plain}

\newtheorem{thm}{Theorem}[section]
\newtheorem{theorem}[thm]{Theorem}
\newtheorem*{theorem*}{Theorem}

\newtheorem*{corollary*}{Corollary}
\newtheorem{lemma}[thm]{Lemma}

\newtheorem*{conjecture*}{Conjecture}

\newtheorem*{question*}{Question}

\newtheorem*{definitions*}{Definitions}
\newtheorem*{convention*}{Convention}
\newtheorem*{conventions*}{Conventions}
\newtheorem*{note*}{Note}
\newtheorem*{exercise*}{Exercise}
\newtheorem*{bibliographical-note*}{Bibliographical note}

\newtheorem*{notation*}{Notation}

\newtheorem*{rem*}{Remark}
\newtheorem{remark}[thm]{Remark}

\newtheorem*{remark*}{Remark}

\newtheorem*{remarks*}{Remarks}
\newtheorem*{example*}{Example}
\newtheorem{example}[thm]{Example}

\newtheorem*{examples*}{Examples}

\newcommand{\A}{\EuScript A}
\newcommand{\barA}{\bar{\A}}
\newcommand{\B}{\EuScript B}
\newcommand{\CC}{\EuScript C}
\newcommand{\DD}{\EuScript D}

\newcommand{\F}{\mathfrak F}
\newcommand{\W}{\mathfrak W}
\newcommand{\GG}{\EuScript G}

\newcommand{\I}{\EuScript I}
\newcommand{\K}{\mathbb K}

\newcommand{\OO}{\EuScript O}
\newcommand{\PP}{\EuScript P}
\newcommand{\QQ}{\EuScript Q}

\renewcommand{\SS}{\EuScript S}
\newcommand{\MM}{\EuScript M}
\newcommand{\NN}{\EuScript N}
\newcommand{\VV}{\EuScript V}
\newcommand{\WW}{\EuScript W}

\newcommand{\ZZ}{\EuScript Z}
\newcommand{\TT}{\EuScript T}
\newcommand{\UU}{\EuScript U}
\renewcommand{\1}{\mathbf{1}}

\newcommand{\boldb}{\mathbf{b}}
\newcommand{\boldo}{\mathbf{o}}
\newcommand{\boldp}{\mathbf{p}}
\newcommand{\boldq}{\mathbf{q}}

\newcommand{\cycl}{\text{\it cycl}}

\newcommand{\semidirect}{\rtimes}
\numberwithin{equation}{section}

\title[$A_\infty$-subalgebras]{$A_\infty$-subalgebras and natural transformations}
\author{Paul Seidel}

\begin{document}
\maketitle

\setcounter{section}{-1}
\section{Introduction}

This paper explores a version of categorical localization. While the results are purely algebraic, the construction itself arose in symplectic topology, specifically in the theory of Lefschetz fibrations, and the connection with localization was arrived at {\em a posteriori}.

The setup is as follows. We consider pairs $(\A,\B)$ consisting of an $A_\infty$-algebra $\B$ and a subalgebra $\A$. We will be primarily working with the category $V = H^0(mod(\A))$ of $\A$-modules. However, that category inherits certain additional data from the presence of $\B$. The most important such data are a functor $F: V \rightarrow V$, given by tensor product with the quotient bimodule $(\B/\A)[-1]$, together with a natural transformation $T: F \rightarrow Id$.

We will introduce a new object associated to $(\A,\B)$, which is an $A_\infty$-algebra $\DD$ with nonvanishing curvature. The definition is quite straightforward, but its meaning is not transparent. We will study an appropriately defined category $W = H^0(modt(\DD))$ of modules over $\DD$. The main result is that $W$ is the localization of $V$ along $T$, a special case of the Verdier quotient construction for triangulated categories. In particular, $W$ only depends on $V$, $F$ and $T$. A similar statement holds for the underlying differential graded categories (see Theorem \ref{th:main} for the precise formulation). We also obtain a corresponding result for the Hochschild homology $HH(\DD)$, which turns out to be expressible in terms of iterated tensor products with $(\B/\A)[-1]$ (Theorem \ref{th:hochschild}).

Objects similar to $\DD$ have already appeared in at least two places in the symplectic literature, namely in the theory of Fukaya categories for closed manifolds \cite{fooo}, and in Chekanov homology (relative Symplectic Field Theory) for Legendrian submanifolds \cite{chekanov99}. We will suggest a tentative link to the latter topic, and further relations to similar symplectic invariants, mostly following \cite{seidel06}.

This manuscript is substantially revised from an earlier version. I would like to thank Tobias Ekholm and Denis Auroux for useful conversations, as well as the referee for a very insightful report. Partial support was provided by NSF grants DMS-0405516 and DMS-0652620.

\section{Localization along a natural transformation\label{sec:loc}}

This section introduces the relevant category theory background. We start in an elementary context, then focus on the triangulated case, and also consider the more modern framework of dg categories. Finally, we give two examples of known constructions which can be interpreted as localizations in our sense.

\subsection*{The basic construction}
Let $V$ be a category. Even though it would be possible to work in complete generality, we prefer to assume (for consistency with later developments) that $V$ is linear over some field $\K$. Let $F: V \rightarrow V$ be an endofunctor, and $T: F \rightarrow Id$ a natural transformation from $F$ to the identity functor. By acting with $F$ on either side, one gets two induced natural transformations $R_FT, L_FT: F^2 \rightarrow F$. Explicitly,
\begin{equation} \label{eq:two-sides}
 (R_FT)_M = T_{F(M)}, \;\;
 (L_FT)_M= F(T_M) \in Hom_V(F^2(M),F(M)).
\end{equation}
These differ in general, but the compositions $T \circ R_FT, \;T \circ L_FT: F^2 \rightarrow Id$ are equal. We denote these compositions by $T^2$, and similarly define higher iterates $T^p: F^p \rightarrow Id$.

We say that $T$ {\em ambidextrous} if the two natural transformations in \eqref{eq:two-sides} agree. Suppose from now on that this is the case. We can then define a new category $Z$, called the {\em localization of $V$ along $T$}, as follows. The objects are the same as in $V$. The morphism spaces are
\begin{equation} \label{eq:formal-limit}
Hom_Z(M,N) = \underrightarrow{\mathit{lim}}_p \, Hom_V(F^p(M),N),
\end{equation}
where the connecting maps in the direct system, $Hom_V(F^p(M),N) \rightarrow Hom_V(F^{p+q}(M),N)$, are given by multiplying with $F^p(T_M^q)$ on the right. Composition of morphisms in $Z$ is inherited from $V$ as follows:
\begin{equation} \label{eq:pq}
\begin{aligned}
& Hom_V(F^q(N),O) \otimes Hom_V(F^p(M),N) \xrightarrow{Id \otimes F^q} \\
& \quad Hom_V(F^q(N),O) \otimes Hom_V(F^{p+q}(M),F^q(N)) \xrightarrow{\text{composition}} Hom_V(F^{p+q}(M),O).
\end{aligned}
\end{equation}
The fact that $T$ is ambidextrous ensures that this is compatible with the direct limit, hence that $Z$ is well-defined.

Our category comes with a functor $I: V \rightarrow Z$, given by the $p = 0$ term in \eqref{eq:formal-limit}. The next result says that this is universal with respect to inverting the morphisms $S = \{T_M\}$, hence that $Z$ agrees with the standard categorical localization $S^{-1}V$.

\begin{lemma}
The image of each $T_M$ under $I$ is an isomorphism. Conversely, each functor from $V$ to another category, which takes the $T_M$ to isomorphisms, factors through $I$ in a unique way.
\end{lemma}

\proof More generally, $T_M^p$ becomes invertible in $Z$, its inverse being the image of the identity morphism $E_{F^p(M)}$ under $Hom_V(F^p(M),F^p(M)) \rightarrow Hom_Z(M,F^p(M))$. An arbitrary morphism $y \in Hom_Z(M,N)$ comes from some $x \in Hom_V(F^p(M),N)$, hence can be written as a composition
\begin{equation} \label{eq:fraction}
\xymatrix{
 M \ar[rr]^-{y} && N \\ &
 \ar[ul]_-{\iso}^-{I(T^p_M)}
 F^p M \ar[ur]_-{I(x)} &
}
\end{equation}
Given that, it is straightforward to factor suitable functors through $Z$. \qed

As a final comment, note that $\{T_M^p\}$ (including the case $p = 0$ which is $T_M^0 = E_M$) is a {\em right localizing class} of morphisms in the sense of \cite{gabriel-zisman}. In that situation, one can define $S^{-1}V$ through a calculus of fractions, which is in fact clearly visible in \eqref{eq:fraction}.

\subsection*{The triangulated case}
Now suppose that $V$ is a triangulated category, and $F$ an exact functor. Let $V_{nil} \subset V$ be the full subcategory consisting of those $M$ such that $T_M^p = 0$ for $p \gg 0$. This is obviously closed under taking cones and direct summands, hence is a thick subcategory \cite[Proposition 1.3]{rickard-stable}. Let $S_{nil}$ be the set of those morphisms whose cone lies in $V_{nil}$. The {\em (Verdier) quotient} category is defined as $V/V_{nil} = S_{nil}^{-1}V$. It carries a canonical induced triangulated structure \cite[p.\ 251]{gelfand-manin}. We now relate this to the previously considered localization $Z = S^{-1}V$.

\begin{lemma} \label{th:ss}
for any $M$, $T_M \in S_{nil}$. Conversely, if $f \in Hom_V(M,N)$ lies in $S_{nil}$, there is a $p$ and morphisms $g \in Hom_V(F^p(N),M)$, $h \in Hom_V(F^p(N),M)$, such that $fg = T_N^p$ and $h F^p(f) = T_M^p$.
\end{lemma}

\proof Taking $C$ to be the mapping cone of $T_M$, we have the commutative diagram
\begin{equation} \label{eq:triangle-d2}
\xymatrix{
 F(M) \ar[r]^-{T_M} & M \ar[r] & C \ar@/_2pc/[ll]_{[1]} \\
 F^2(M) \ar[r]_-{F(T_M)}  \ar[u]^-{T_{F(M)}} &
 F(M) \ar[u]^-{T_M} \ar[r] \ar@{-->}[ur] &
 F(C) \ar[u]_-{T_C} \ar@/^2pc/[ll]^{[1]}
}
\end{equation}
The dashed arrow is obviously zero, hence $T_C$ factors through $M$. After applying $F$ one sees that $T_C^2 = 0$.

For the converse statement, take $f$, and let $C$ be its cone. Look at the long exact sequence
\begin{equation}
 \cdots \rightarrow Hom_V(F^p(N),M) \xrightarrow{f \cdot} Hom_V(F^p(N),N) \longrightarrow Hom_V(F^p(N),C) \rightarrow \cdots
\end{equation}
The image of $T_N^p$ under the second map agrees with that of $T_C^p$ under $Hom_V(F^p(C),C) \rightarrow Hom_V(F^p(N),C)$, hence vanishes for large $p$. Taking a preimage under the first map produces the desired $g$. The same kind of argument applies to multiplication by $F^p(f)$ on the right. \qed

Lemma \ref{th:ss} shows that $S \subset S_{nil}$, and conversely that any element of $S_{nil}$ becomes invertible once we invert $S$ (to get a left inverse, take the equality $h F^p(f) = T_M^p$ and write this as $h (T_N^p)^{-1} T_N^p F^p(f) = h (T_N^p)^{-1} f T_M^p = T_M^p$). We therefore get functors in both directions, which yield an equivalence
\begin{equation} \label{eq:compare-quotient}
Z = S^{-1}V \iso S_{nil}^{-1}V = V/V_{nil}.
\end{equation}

\subsection*{The differential graded case}
Next, assume that $V = H^0(\VV)$ is the cohomology category associated to a pre-triangulated dg (differential graded) category. Given $F$ and $T$ as before, one denotes by $\VV_{nil} \subset \VV$ the full dg subcategory with $H^0(\VV_{nil}) = V_{nil}$. In this setup, one can form the {\em dg quotient} category $\ZZ = \VV/\VV_{nil}$ in the sense of \cite{drinfeld02} (simplifying an earlier construction of \cite{keller99a}). On the cohomology level, $H^0(\ZZ) \iso V/V_{nil}$ recovers our previous $Z$, by \cite[Theorem 3.4]{drinfeld02}. It also comes with a natural dg functor $\I: \VV \rightarrow \ZZ$ such that $H^0(\I) = I$. As a consequence of our previous discussion, there is an easy criterion for identifying when another given dg category is quasi-equivalent to this quotient:

\begin{lemma} \label{th:dg-quotient}
Let $\WW$ be a pre-triangulated dg category, and $\GG: \VV \rightarrow \WW$ a dg functor. Write $W$, $G$ for the corresponding cohomology level structures. Suppose that: (i) the objects in the image of $G$ generate $W$ as a triangulated category; (ii) $G(T_M)$ is invertible for all $M$; (iii) the resulting maps
\begin{equation} \label{eq:lim-map}
\underrightarrow{\mathit{lim}}_p \, Hom_V(F^p(M),N) \longrightarrow Hom_W(G(M),G(N))
\end{equation}
are isomorphisms. In that case, $\WW$ is quasi-equivalent to $\ZZ = \VV/\VV_{nil}$.
\end{lemma}

\proof By \eqref{eq:lim-map}, every $x \in Hom_W(G(M),G(N))$ can be written as a composition
\begin{equation}
G(M) \xrightarrow{G(T_M^p)^{-1}} G(F^p(M)) \xrightarrow{G(y)} G(M)
\end{equation}
for some $y$, where the first arrow is an isomorphism. Note that because it comes from an underlying dg functor, $G$ is automatically exact. It follows that the image of the cone of $y$ under $G$ is isomorphic to the cone of $x$. This proves that the objects in the image of $G$ form a triangulated subcategory, which together with the other assumption shows that $G$ is essentially surjective (onto on isomorphism classes of objects).

$G$ inverts the $T_M$, hence by Lemma \ref{th:ss} all of $S_{nil}$, which means that it kills all of $V_{nil}$. By definition of dg quotient, we therefore have an induced dg functor $\ZZ \rightarrow \WW$. By looking at \eqref{eq:lim-map} and \eqref{eq:compare-quotient}, one sees that this functor is fully faithful. On the other hand, we have shown that the underlying cohomology level functor is essentially surjective, so we have a quasi-equivalence. \qed

\subsection*{Complements of divisors}
Let $X$ be a smooth quasiprojective variety over $\K = \C$, $Y \subset X$ a (possibly singular) hypersurface, and $U = X \setminus Y$ its complement. We will be interested in the relation between coherent sheaves on $X$ and on $U$. Write $Y = s^{-1}(0)$, where $s$ is the canonical section of the line bundle $L = {\EuScript O}_X(Y)$. The basic results are as follows:
\begin{equation} \label{eq:poles}
\parbox{36em}{
For any coherent sheaf $M$ on $X$, $H^*(M|U) \iso \underrightarrow{\mathit{lim}}_p H^*(M \otimes L^p)$, where the direct limit is formed with respect to multiplication with $s$.
}
\end{equation}
This is obvious for $H^0$, and for higher $H^k$ it can be proved by looking at the Cech complexes associated to affine covers.
\begin{equation} \label{eq:extend}
\parbox{36em}{
Every coherent sheaf on $U$ is isomorphic to the restriction of a coherent sheaf on $X$.
}
\end{equation}
To see that, note first that every coherent sheaf on $U$ can be written as a quotient of a map of vector bundles. More precisely, one can take both vector bundles to be direct sums of powers of some fixed ample line bundle. Taking that line bundle to be the restriction of an ample line bundle on $X$, we have written our sheaf as the quotient of a map $M|U \rightarrow N|U$, where $M,N$ live on $X$. That map can have poles along $Y$, but one can get rid of those by replacing $N$ with $N \otimes L^p$ for $p \gg 0$.

We now pass to the level of bounded derived categories of coherent sheaves, where the basic results are:
\begin{align} \label{eq:derived-extend}
\parbox{36em}{
The restriction functor $D^b(X) \rightarrow D^b(U)$ is essentially surjective.} \\ \label{eq:ext-lim}
\parbox{36em}{
For any $M,N \in Ob\,D^b(X)$, $Hom_{D^b(U)}(M|U,N|U) \iso \underrightarrow{\mathit{lim}}_p \, Hom_{D^b(X)}(M,N \otimes L^p)$.}
\end{align}
The first statement is obvious from \eqref{eq:extend} and \eqref{eq:poles} above. The second one reduces to \eqref{eq:poles} in the case where $M$ and $N$ are vector bundles or shifts thereof, via $Ext^i(M,N) = H^i(M^\vee \otimes N)$. Any object of the bounded derived category can be represented by a bounded complex of vector bundles, so a filtration argument establishes the general case (we point out that \eqref{eq:ext-lim} fails for quasi-coherent sheaves, such as $\bigoplus_{k =0 }^{\infty} \OO_X$, as well as for unbounded complexes of coherent sheaves, such as $\bigoplus_{k = 0}^\infty \OO_X[\pm k]$).

To fit this into our general context, take $V = D^b(X)$, $F$ the functor of tensoring with $L^{-1}$, and $T: F \rightarrow Id$ the natural transformation given by multiplying with $s$. We then have an induced exact functor from the localized category $Z$ to $W = D^b(U)$, which is an equivalence by \eqref{eq:derived-extend} and \eqref{eq:ext-lim}. In this case, the thick subcategory $V_{nil}$ consists precisely of those complexes whose cohomology is supported on $Y$ (these are actually isomorphic to complexes living on a formal thickening of $Y$, see \cite[Lemma 3]{bezrukavnikov00}).

\subsection*{Landau-Ginzburg branes}
We consider the same situation as before, but concentrate on $Y$. Recall that an object of $D^b(Y)$ is called a {\em perfect complex} if it is isomorphic to a bounded complex of vector bundles.
There is a useful cohomological criterion \cite[Lemma 1.11]{orlov-equi}:
\begin{equation} \label{eq:cohomology-criterion}
\parbox{36em}{$M$ is perfect if and only if for all $N \in Ob\,D^b(Y)$, $Hom_{D^b(Y)}(M,N[d]) = 0$ provided that $d \gg 0$.}
\end{equation}
This in particular shows that perfect complexes form a thick triangulated subcategory of $D^b(Y)$.
The quotient by that subcategory is called the {\em category of Landau-Ginzburg branes} on $Y$, denoted by $D^b_{sing}(Y)$ \cite{orlov-branes}.

To interpret this in our framework, let $j: Y \rightarrow X$ be the inclusion, and $j_*, j^*$ the (derived) pushforward and pullback functors. For any $M \in Ob\,D^b(Y)$ there is a canonical distinguished triangle
\begin{equation} \label{eq:push-pull}
\xymatrix{
 M \otimes (L^{-1}|Y)[1] \ar[r] & j^*j_*M \ar[r] & M \ar@/^2pc/[ll]^{}
}
\end{equation}
There is a more universal viewpoint on this, as follows. In $D^b(Y \times Y)$ we have a triangle
\begin{equation} \label{eq:fm-sequence}
\xymatrix{
 {\EuScript O}_{\Delta_Y} \otimes (L^{-1}|Y)[1] \ar[r] &
 (j \times j)^* {\EuScript O}_{\Delta_X} \ar[r] &
 {\EuScript O}_{\Delta_Y}
 \ar@/^2pc/[ll]^{}
}
\end{equation}
where $\Delta$ stands for the diagonals. The induced action on $D^b(Y)$ by Fourier-Mukai transforms yields \eqref{eq:push-pull}. Set $V = D^b(Y)$, and let $F: V \rightarrow V$ be the exact functor of tensoring with $(L|Y)[-2]$. The boundary maps of \eqref{eq:push-pull}, inherited from the corresponding one in \eqref{eq:fm-sequence}, constitute an ambidextrous natural transformation $T: F \rightarrow Id$, which we can use to build the localized category $Z$.

To relate the two constructions, we need the following observation:
\begin{equation} \label{eq:nil}
\parbox{36em}{$M$ is perfect iff the $p$-th iterate of the boundary map of \eqref{eq:push-pull}, which is an element of $Hom_{D^b(Y)}(M,M \otimes L^{-p}(Y)[2p])$, vanishes for $p \gg 0$.}
\end{equation}
Namely, if $M$ is perfect and $r \gg 0$, then $Hom_{D^b(Y)}(M, M \otimes (L^{-p}|Y)[r]) = 0$ for all $p$, so the $p$-th iterate will vanish for degree reasons. In converse direction, we know that the cone of the boundary map is a perfect complex, and the same then holds for its iterates. If such an iterate is zero, $M$ itself is a direct summand of a perfect complex, hence again perfect. The upshot is that the category $V_{nil}$ is exactly that of perfect complexes, hence that $Z$ is equivalent to $D^b_{sing}(Y)$. As a consequence,
\begin{equation}
Hom_{D^b_{sing}(Y)}(M,N) = \underrightarrow{\mathit{lim}}_p\, Hom_{D^b(Y)}(M, N
\otimes (L^{-p}|Y)[2p]).
\end{equation}

\section{$A_\infty$-bimodules}
The notion of bimodule over an $A_\infty$-algebra appears in many places in the literature, see for instance \cite{tradler01,ks-notes,ukraine}. Nevertheless, it may still be appropriate to give a concise presentation, if only to clarify the conventions.

\subsection*{Notation}
As before, everything will be linear over a fixed ground field $\K$. An $A_\infty$-algebra is a graded vector space $\A$ together with a sequence of multilinear maps $\mu_\A^d: \A^{\otimes d} \rightarrow \A$, $d \geq 1$, satisfying certain associativity equations. More compactly, if we set $T(\A[1]) = \K \oplus \A[1] \oplus \A[1]^{\otimes 2} \oplus \cdots$, the $A_\infty$-structure can be written as a single graded map
\begin{equation} \label{eq:ainfty}
\mu_\A: T(\A[1]) \longrightarrow \A[2],
\end{equation}
with vanishing constant term $\K \rightarrow \A[2]$. We will use the standard ``reduced'' or ``bar'' sign convention, so the $A_\infty$-associativity equations are
\begin{equation} \label{eq:a}
\sum_{m,n} (-1)^\ast \mu_\A^{s-n+1}(a_1,\dots,a_m,\mu_\A^n(a_{m+1},\dots,a_{m+n}),a_{m+n+1},\dots,a_s) = 0,
\end{equation}
where $\ast = \|a_1\|+\cdots+\|a_m\|$ is the sum of reduced degrees $\|a_i\| = |a_i|-1$. To memorize the sign, think of $\mu_\A$ as acting from the left, so that $\mu_\A^n$ (which is of degree $1$ in the reduced grading) has to be commuted past $a_1,\dots,a_m$. From now on, {\em given any expression which consists of a multilinear map applied to a block of subsequent entries in a tensor expression, we write $\ast$ for the sum of degrees of all entries lying to the left of that block, where ``degree'' for elements of $A_\infty$-algebras means reduced degree}.
Further, all our $A_\infty$-algebras are assumed to be strictly unital, which means that there is a (necessarily unique) $e \in \A^0$ such that
\begin{equation}
\mu^1_\A(e) = 0, \;\; \mu^2_\A(e,a) = a, \;\; \mu^2_\A(a,e) = (-1)^{\|a\|+1} a, \;\;
\mu^d_\A(\dots,e,\dots) = 0 \text{ for any $d>2$.}
\end{equation}
We write $\bar{\A} = \A/\K e$.

\subsection*{Bimodules}
Let $\A,\A'$ be two $A_\infty$-algebras. An $(\A,\A')$-bimodule consists of a graded vector space $\PP$ together with a structure map
\begin{equation} \label{eq:bimodule-structure}
\mu_\PP: T(\A[1]) \otimes \PP \otimes T(\A'[1]) \longrightarrow \PP[1],
\end{equation}
satisfying a generalized bimodule equation. Write $\mu_\PP^{r|1|s}$ for the component of \eqref{eq:bimodule-structure} which has $r$ tensor factors on the left and $s$ on the right. Then, the $A_\infty$-bimodule equation is
\begin{equation} \label{eq:bimodule}
\begin{aligned}
 & \sum_{m,n}
 (-1)^\ast
 \mu_\PP^{r|1|s-n+1}(a_1,\dots,a_r,\boldp,a'_1,\dots,a'_m,\mu_{\A'}^n(a'_{m+1},\dots,a'_{m+n}), a'_{m+n+1},\dots,a'_s)
 \\ + \;
 & \sum_{m,n} (-1)^\ast
 \mu_\PP^{r-n+1|1|s}(a_1,\dots,a_m,\mu_\A^n(a_{m+1},\dots,a_{m+n}),a_{m+n+1},\dots,a_r,\boldp,a'_1,\dots,
 a_s')
 \\ + \;
 & \sum_{m,n} (-1)^\ast
 \mu_\PP^{m|1|s-n}(a_1,\dots,a_m,\mu_\PP^{r-m|1|n}(a_{m+1},\dots,\boldp,\dots,a'_n),a'_{n+1},\dots,a'_s)
 \; = \; 0.
\end{aligned}
\end{equation}
The boldface notation for elements of $\PP$ is just an attempt to increase readability in long formulae. Also, we follow the previous indications for signs $(-1)^\ast$, with the proviso that the degree of $\boldp$ is the natural (unreduced) one. Finally, we impose a unitality requirement, which is
\begin{equation}
\begin{aligned}
 & \mu_\PP^{1|1|0}(e,\boldp) = \boldp, \;\; \mu_\PP^{0|1|1}(\boldp,e') = (-1)^{|\boldp|+1} \boldp, \\
 & \mu_\PP^{r|1|s}(\dots,e,\dots,\boldp,\dots) =
 \mu_\PP^{r|1|s}(\dots,\boldp,\dots,e',\dots) =
 0 \text{ for $r+s>1$.}
\end{aligned}
\end{equation}

\subsection*{Morphisms of bimodules}
Let $\PP$, $\QQ$ be two $(\A,\A')$-bimodules. A pre-homomorphism $\rho: \PP \rightarrow \QQ$ of degree $|\rho| = k$ is a map
\begin{equation} \label{eq:pre-hom}
\begin{aligned}
&
\rho: T(\A[1]) \otimes \PP \otimes T(\A'[1]) \longrightarrow \QQ[k], \\
& \rho^{r|1|s}(\dots,e,\dots,\boldp,\dots) = \rho^{r|1|s}(\dots,\boldp,\dots,e',\dots) = 0.
\end{aligned}
\end{equation}
The unitality condition in \eqref{eq:pre-hom} can equivalently be expressed by writing $\rho$ as a map $T(\bar{\A}[1]) \otimes \PP \otimes T(\bar{\A}'[1]) \rightarrow \QQ[k]$. We call $\rho$ a bimodule homomorphism if it satisfies $\partial\rho = 0$, where
\begin{equation} \label{eq:bimodule-d}
\begin{aligned}
 & (\partial\rho)^{r|1|s}(a_1,\dots,a_r,\boldp,a'_1,\dots,a'_s) =
 \\ &
 \sum_{m,n}
 (-1)^{|\rho| \ast} \
 \mu_\QQ^{m|1|s-n}(a_1,\dots,a_m,\rho^{r-m|1|n}(a_{m+1},\dots,\boldp,\dots,a'_n),a'_{n+1},\dots,a'_s)
 \\ + &
 \sum_{m,n}
 (-1)^{|\rho|+1+\ast}
 \rho^{m|1|s-n}(a_1,\dots,a_m,\mu_\PP^{r-m|1|n}(a_{m+1},\dots,\boldp,\dots,a'_n),a'_{n+1},\dots,a'_s)
 \\ + &
 \sum_{m,n}
 (-1)^{|\rho|+1+\ast}
 \rho^{r|1|s-n+1}(a_1,\dots,\boldp,\dots,a_m',\mu_{\A'}^n(a'_{m+1},\dots,a'_{m+n}),a'_{m+n+1},\dots,a'_s)
 \\ + &
 \sum_{m,n}
 (-1)^{|\rho|+1+\ast}
 \rho^{r-n+1|1|s}(a_1,\dots,a_m,\mu_\A^n(a_{m+1},\dots,a_{m+n}),a_{m+n+1},\dots,\boldp,\dots,a'_s).
\end{aligned}
\end{equation}
Composition is defined by a similar, if somewhat simpler, formula:
\begin{equation} \label{eq:bimodule-c}
\begin{aligned}
 & (\theta\rho)^{r|1|s}(a_1,\dots,a_r,\boldp,a'_1,\dots,a'_s) = \\ &
 \sum_{m,n} (-1)^{|\rho|\ast}
 \theta^{n|1|s-m}(a_1,\dots,a_m,\rho^{r-m|1|n}(a_{m+1},\dots,\boldp,\dots,a'_n),a'_{n+1},\dots,a'_s).
\end{aligned}
\end{equation}
A homomorphism whose higher order terms all vanish will be called naive. In the case of degree zero, this consists of a single graded linear map $\rho^{0|1|0}$ which strictly commutes with all bimodule operations; and such homomorphisms are composed in the obvious way. A trivial example is the identity endomorphisms $\1 = \1_{\PP}$, given by $\1^{0|1|0}(\boldp) = \boldp$, with vanishing higher order terms.

There is a dg category $\UU = mod(\A,\A')$ whose objects are bimodules, and whose morphisms are pre-homomorphisms. Any bimodule homomorphism $\rho: \PP \rightarrow \QQ$ induces a map $H(\rho): H(\PP) \rightarrow H(\QQ)$, where the cohomology is taken with respect to $\mu^{0|1|0}$. If this induced map is an isomorphism, $\rho$ is called a quasi-isomorphism. An important property of the $A_\infty$ theory is that any quasi-isomorphism is an isomorphism in $H(\UU)$. Concretely, this means that there is another bimodule homomorphism $\tau: \QQ \rightarrow \PP$ such that $\tau\rho$ and $\rho\tau$ are homotopic (in the dg category sense) to the respective identities. The proof is a standard argument involving the filtration of $hom_\UU$ spaces by length.

\subsection*{Shift}
There is a natural shift operation on bimodules, $\PP[1]^i = \PP^{i+1}$. This is accompanied by a sign change in the bimodule structure maps. More precisely,
\begin{equation} \label{eq:p-shift}
\mu_{\PP[1]}^{r|1|s}(a_1,\dots,a_r,\boldp,a'_1,\dots,a'_s) = (-1)^{\circ+1}
\mu_{\PP}^{r|1|s}(a_1,\dots,a_r,\boldp,a'_1,\dots,a'_s).
\end{equation}
Here $\circ = \|a_1\| + \cdots + \|a_r\|$. As in the case of $\ast$, we will promote this to standard notation: {\em whenever we are dealing with some particular bimodule operation, $\circ$ is the sum of reduced degrees for the elements of $\A$ lying to the left of the $\PP$ factor.} As a consequence of the sign convention in \eqref{eq:p-shift}, the naive identification $hom_\UU^{i+1}(\PP,\QQ) = hom_\UU^i(\PP,\QQ[1])$ is an isomorphism of chain complexes. In contrast, the isomorphism $hom_\UU^{i-1}(\PP,\QQ) \iso hom_\UU^i(\PP[1],\QQ)$, which we write as $\rho \mapsto \rho[1]$, should be defined as $(\rho[1])^{r|1|s} = (-1)^\circ \rho^{r|1|s}$.

\subsection*{Short exact sequences}
Let $\QQ$ be an $(\A,\A')$-bimodule, and $\PP \subset \QQ$ a submodule, which means a graded subspace closed under all the operations $\mu^{r,1,s}_\QQ$. This obviously inherits a bimodule
structure, and so does the quotient $\OO = \QQ/\PP$. The inclusion and projection maps $\iota: \PP \rightarrow \QQ$, $\pi: \QQ \rightarrow \OO$, can be promoted to bimodule homomorphisms in the naive way. We call the outcome a {\em naive short exact sequence} of bimodules,
\begin{equation} \label{eq:short}
 0 \rightarrow \PP \stackrel{\iota}{\longrightarrow} \QQ \stackrel{\pi}{\longrightarrow} \OO \rightarrow
 0.
\end{equation}
Following the standard pattern, one can construct a boundary operator $\delta$, which is a homomorphism $\OO \rightarrow \PP$ of degree one. To define this, choose a splitting $\sigma$ of $\pi$ as a map of graded vector spaces. Extend that naively to a pre-homomorphism of bimodules, which means an element of $hom_\UU^0(\OO,\QQ)$. Then take $\delta = \partial\sigma$, which will automatically lie in the subspace $hom_\UU^1(\OO,\PP) \subset hom_\UU^1(\OO,\QQ)$. The explicit formula is
\begin{equation} \label{eq:boundary-map}
 \delta^{r|1|s}(a_1,\dots,a_r,\boldo,a'_1,\dots,a'_s) = (id-\sigma\pi)
 \mu_\QQ^{r|1|s}(a_1,\dots,a_r,\sigma(\boldo),a_1',\dots,a_s').
\end{equation}
Note that $(id-\sigma\pi)$ is the induced splitting of $\iota$. Clearly, the homotopy class of $\delta$ is independent of the choice of splitting.

In converse direction, given a unital degree $1$ bimodule homomorphism $\delta: \OO \rightarrow \PP$, we define its mapping cone to be $\QQ = \OO \oplus \PP$, with the bimodule structure given by
\begin{equation}
 \mu_\QQ^{r|1|s} = \begin{pmatrix} \mu_\OO^{r|1|s} & 0 \\ \delta^{r|1|s} & \mu_\PP^{r|1|s}
 \end{pmatrix}.
\end{equation}
This fits into an obvious short exact sequence, whose boundary map is the given homomorphism $\delta$.

\subsection*{Tensor product}
Let $\PP$ be an $(\A,\A')$-bimodule, and $\QQ$ an $(\A',\A'')$-bimodule. We define their tensor product to be the $(\A,\A'')$-bimodule
\begin{equation} \label{eq:tensor}
 \SS = \PP \otimes_{\A'} \QQ \stackrel{\mathrm{def}}{=} \PP \otimes T(\barA'[1]) \otimes \QQ,
\end{equation}
with the differential
\begin{align} \label{eq:010-tensor}
& \begin{aligned}
 & \mu_\SS^{0|1|0}(\boldp \otimes a'_1 \otimes \cdots \otimes a'_t \otimes \boldq)
 = \\
 & \sum_n \mu_\PP^{0|1|n}(\boldp,a'_1,\dots,a'_n) \otimes a'_{n+1} \otimes \cdots
 \otimes a'_t \otimes \boldq \\
 + & \sum_n (-1)^\ast \boldp \otimes a'_1 \otimes \cdots \otimes a'_n \otimes
 \mu_\QQ^{t-n|1|0}(a'_{n+1},\dots,a'_t,\boldq) \\
 + & \sum_{m,n} (-1)^\ast \boldp \otimes a'_1 \otimes \cdots \otimes a'_m
 \otimes \mu_{\A'}^n(a'_{m+1},\dots,a'_{m+n}) \otimes a'_{m+n+1} \otimes \cdots \otimes a'_t \otimes \boldq,
\end{aligned}
\intertext{the one-sided higher order structure maps} &
\begin{aligned}
 & \mu_\SS^{r|1|0}(a_1,\dots,a_r,\boldp \otimes a'_1 \otimes \cdots \otimes a'_t \otimes \boldq) = \\
 & \quad \sum_n \mu_\PP^{r|1|n}(a_1,\dots,a_r,\boldp,a'_1,\dots,a'_n) \otimes a'_{n+1} \otimes \cdots
 \otimes a'_t \otimes \boldq, \\
 & \mu_\SS^{0|1|s}(\boldp \otimes a_1' \otimes \cdots \otimes a'_t \otimes \boldq, a''_1,\dots,a''_s) = \\
 & \quad
 \sum_n (-1)^\ast \boldp \otimes a'_1 \otimes \cdots \otimes a'_n \otimes \mu_\QQ^{t-n|1|s}(
 a'_{n+1},\dots,a'_t,\boldq,a''_1,\dots,a''_s);
\end{aligned}
\intertext{and finally,} &
 \mu_\SS^{r|1|s} = 0 \quad \text{whenever $r>0$ and $s>0$.}
 \label{eq:tensor-vanish}
\end{align}
In the classical case (for complexes of bimodules over associative algebras), the formulae \eqref{eq:tensor} specialize to the derived tensor product, with the middle term $T(\bar\A'[1])$ arising from the reduced bar resolution of the diagonal bimodule. As a philosophical aside, note that in the $A_\infty$ world any notion of tensor product must necessarily be a ``derived'' one, because of the invertibility of quasi-isomorphisms.

Tensor products are functorial in the obvious sense. If $\tau: \PP \rightarrow \PP'$ and $\rho: \QQ \rightarrow \QQ'$ are bimodule pre-homomorphisms, one defines $\tau \otimes \rho: \PP \otimes_{\A'} \QQ \rightarrow \PP' \otimes_{\A'} \QQ'$ by
\begin{equation}
\begin{aligned}
 & (\tau \otimes \rho)^{r|1|s}(a_1,\dots,a_r,\boldp \otimes a'_1 \otimes \cdots \otimes a'_t \otimes
 \boldq,a_1'',\dots,a_s'') \\ & =
 \sum_{m,n}  (-1)^{|\rho| \ast}
 \tau^{r|1|m}(a_1,\dots,\boldp,\dots,a'_m) \otimes a'_{m+1} \otimes \cdots \otimes a'_n
 \otimes \rho^{t-n|1|s}(a'_{n+1},\dots,\boldq,\dots,a''_s).
\end{aligned}
\end{equation}
To prevent confusion, we spell out what the standard sign means in this case, namely $\ast = \|a_1\| + \cdots + \|a_r\| + |\boldp| + \|a'_1\| + \cdots + \|a'_n\|$. This construction has the expected properties
\begin{equation}
\begin{aligned}
 & \partial(\tau \otimes \rho) = \partial\tau \otimes \rho + (-1)^{|\tau|} \tau \otimes
 \partial\rho, \\
 & (\tau_1 \otimes \rho_1) (\tau_2 \otimes \rho_2) = (-1)^{|\rho_1|\,|\tau_2|}
 \tau_1\tau_2 \otimes \rho_1\rho_2.
\end{aligned}
\end{equation}

\subsection*{Tensor product with the diagonal}
Following standard usage, we talk of $\A$-bimodules instead of $(\A,\A)$-bimodules. The standard example is the diagonal bimodule, whose underlying graded vector space is $\A$, with operations
\begin{equation} \label{eq:diagonal}
\mu^{r|1|s}_\A = (-1)^{\circ+1} \mu_\A^{r+1+s}.
\end{equation}
The shifted version $\A[1]$ is actually a little simpler, since $\mu_{\A[1]}^{r|1|s} = \mu_\A^{r+1+s}$.

Naively, one expects the tensor product with the diagonal to do nothing. In the present framework, this holds only up to quasi-isomorphism. More precisely, we have a natural quasi-isomorphism
\begin{equation} \label{eq:epsilon-left}
\begin{aligned}
 & \epsilon = \epsilon_{\PP,left}: \A \otimes_\A \PP \longrightarrow \PP, \\
 & \epsilon^{r|1|s} = (-1)^\circ \mu_\PP^{r+1+t|1|s}: \barA[1]^{\otimes r} \otimes \A \otimes \barA[1]^{\otimes t} \otimes \PP \otimes \barA[1]^{\otimes s} \longrightarrow \PP
\end{aligned}
\end{equation}
(as before, if we shift to $\epsilon[1]: (\A \otimes_\A \PP)[1] = \A[1] \otimes_\A \PP \rightarrow \PP$, the sign vanishes). Naturality of $\epsilon$ means that for any bimodule homomorphism $\rho: \PP \rightarrow \QQ$, one has a homotopy commutative diagram
\begin{equation} \label{eq:a-fun}
\xymatrix{
 \A \otimes_\A \PP \ar[rr]^-{\epsilon_{\PP,left}} \ar[d]^{\1_\A \otimes \rho} && \PP \ar[d]^{\rho} \\
 \A \otimes_\A \QQ \ar[rr]^-{\epsilon_{\QQ,left}} && \QQ
}
\end{equation}
(in the notation from \eqref{eq:epsilon-left}, the homotopy is $\kappa^{r|1|s} = (-1)^{\circ} \rho^{r+1+t|1|s}$). The fact that $\epsilon$ is a quasi-isomorphism follows from standard spectral sequence arguments, which reduce things to the cohomological level and hence to the classical case of algebras. One can also construct an explicit quasi-inverse $\xi = \xi_{\PP,left}: \PP \rightarrow \A \otimes_\A \PP$. Namely, take
\begin{equation} \label{eq:quasi-inverse}
\xi^{r|1|0}(a_1,\dots,a_r,\boldp) = e \otimes a_1 \cdots \otimes a_r \otimes \boldp,
\end{equation}
and set the other terms $\xi^{r|1|s}$, $s>0$, to zero. This is actually strictly natural, meaning that the diagram corresponding to \eqref{eq:a-fun} (with the horizontal arrows pointing left) commutes on the nose.

For completeness, we should also mention the right-sided counterparts of \eqref{eq:epsilon-left} and \eqref{eq:quasi-inverse},
\begin{equation} \label{eq:epsilon-right}
\begin{aligned}
 & \epsilon = \epsilon_{\PP,right}: \PP \otimes_{\A' } \A' \longrightarrow \PP, \\
 & \xi = \xi_{\PP,right}: \PP \longrightarrow \PP \otimes_{\A'} \A'.
\end{aligned}
\end{equation}
The first is given by $\epsilon^{r|1|s} = (-1)^\ddag \mu_\PP^{r|1|s+1+t}$, where $\ddag = \|a_1\| + \cdots + \|a_r\| + |\boldp| + \|a'_1\| + \cdots + \|a'_t\| + 1$. Similarly, its quasi-inverse has nontrivial components $\xi^{0|1|s}(\boldp,a'_1,\dots,a'_s) = \boldp \otimes a'_1 \otimes a'_s \otimes e'$. Given any $(\A,\A')$-bimodule $\PP$ and $(\A',\A'')$-bimodule $\QQ$, we have two possible choices of quasi-isomorphisms in either direction:
\begin{equation} \label{eq:left-right-choice}
\begin{aligned}
& \1_{\PP} \otimes \epsilon_{\QQ,left},\, \epsilon_{\PP,right} \otimes \1_{\QQ}:
\PP \otimes_{\A'} \A' \otimes_{\A'} \QQ \longrightarrow \PP \otimes_{\A'} \QQ, \\
& \1_{\PP} \otimes \xi_{\QQ,left}, \, \xi_{\PP,right} \otimes \1_{\QQ}: \PP \otimes_{\A'} \QQ \longrightarrow \PP \otimes_{\A'} \A' \otimes_{\A'} \QQ.
\end{aligned}
\end{equation}
The first pair turn out to be homotopic, while the second pair agree strictly. Note also that in the special case of the diagonal module, the left-sided and right-sided $\epsilon$ maps give rise to the same homomorphism $\A \otimes_\A \A \rightarrow \A$, while the quasi-inverses only yield homotopic homomorphisms $\A \rightarrow \A \otimes_\A \A$.

\subsection*{Convolution functors}
An $A_\infty$-module over $\A$ is the same as an $(\A,\K)$-bimodule. In this situation, we will specialize the notation to make it more familiar, writing $\mu^{s|1}$ rather than $\mu^{s|1|0}$, and $mod(\A)$ rather than $mod(\A,\K)$. These are of course left modules; right $\A$-modules, defined as objects of $mod(\K,\A)$, will also arise occasionally.

Take two $A_\infty$-algebras $\A$ and $\A'$, and consider the associated dg categories $\VV = mod(\A)$, $\VV' = mod(\A')$. Every $\PP \in Ob\,mod(\A,\A')$ gives rise to a dg functor $\Phi_\PP: \VV' \rightarrow \VV$, which sends a module $\MM$ to $\PP \otimes_{\A'} \MM$, and a module pre-homomorphism $\phi$ to $\1_\PP \otimes \phi$. We call this {\em convolution with $\PP$}. Tensor
product of bimodules corresponds to composition of convolution functors. Moreover, a bimodule homomorphism $\rho: \PP \rightarrow \QQ$ gives rise to a natural transformation $\Phi_\rho: \Phi_\PP \rightarrow \Phi_\QQ$, which consists of the collection of module homomorphisms
$\rho \otimes \1_\MM$. A more precise formulation would be as follows: let $fun(\VV',\VV)$ be the dg category of all dg functors and their natural (pre-)transformations. Then convolution defines a canonical dg-functor
\begin{equation} \label{eq:functorial-phi}
\Phi: mod(\A,\A') \longrightarrow fun(\VV',\VV).
\end{equation}
There is one rather unfortunate aspect about this formalism. Taking $\PP = \A$ to be the diagonal bimodule, one finds that the maps \eqref{eq:quasi-inverse} provide a natural transformation $Id \longrightarrow \Phi_{\A}$. By definition, each single element making up that natural transformation is an isomorphism in $H^0(\VV)$. However, it is not clear that the natural transformation itself is an isomorphism in $H^0(fun(\VV,\VV))$, since the inverse maps \eqref{eq:epsilon-left} are not functorial on the cochain level. One can improve the situation by considering a suitably modified functor category (either adjoining abstract formal inverses \`a la derived category, or more concretely using $A_\infty$-natural transformations). We will usually adopt a more crude solution, which is to stay on the level of $H^0(\VV)$, where the problem is of course non-existent.

\subsection*{Bar constructions}
Assume that $\A$ is augmented, which means that it can be written as
\begin{equation}
\A = \K e \oplus \bar{\A},
\end{equation}
with $\bar{\A}$ itself a non-unital $A_\infty$-algebra. Consider $T(\barA[1])$ as a free graded coalgebra, with the standard coproduct. $\mu_\A$ determines a coalgebra differential of this, given explicitly by
\begin{equation} \label{eq:bar-d}
\partial_{T(\barA[1])}(a_1 \otimes \cdots \otimes a_s) = \sum_{m,n} (-1)^\ast a_1 \otimes \cdots \otimes \mu_{\A}^n(a_{m+1},\dots,a_{m+n}) \otimes a_{m+n+1} \otimes \cdots \otimes a_s.
\end{equation}
Similarly, given two augmented $A_\infty$-algebras $\A$ and $\A'$ and a bimodule $\PP$, the bimodule structure gives rise to a differential on $T(\barA[1]) \otimes \PP \otimes T(\barA'[1])$, which makes that space into a bi-comodule over the associated dg coalgebras. This construction in fact yields a full embedding of the dg category of $A_\infty$-bimodules into that of dg bi-comodules, which is also compatible with tensor products. However, if we start with the diagonal $A_\infty$-bimodule, then $T(\barA[1]) \otimes \A \otimes T(\barA[1])$ is not the diagonal dg bimodule, but only a resolution thereof. This is one way of explaining the difficulties encountered above when tensoring with the diagonal.

%

\section{$A_\infty$-subalgebras\label{sec:subalgebra}}

This section is the computational core of the paper. We introduce one of the two key objects, namely the boundary homomorphism associated to a pair $\A \subset \B$ of $A_\infty$-algebras. We also carry out some related computations with $\A$-modules and $\A$-bimodules. These reveal cancellation phenomena which will be crucial later on.

\subsection*{The boundary homomorphism}
From now on, we work in the following context: $\B$ is an $A_\infty$-algebra, and $\A$ is an $A_\infty$-subalgebra (containing the identity element of $\B$). Consider the naive short exact sequence of $\A$-bimodules
\begin{equation} \label{eq:b-a-sequence}
0 \rightarrow \A \stackrel{\iota}{\longrightarrow} \B \stackrel{\pi}{\longrightarrow} \B/\A \rightarrow 0.
\end{equation}
Here $\A$ is the diagonal bimodule, and $\B$ is similarly considered to be an $\A$-bimodule by
restriction of the diagonal $\B$-bimodule structure. Let $\delta$ be the boundary homomorphism of this
exact sequence. Explicitly, choosing a splitting $\sigma: \B/\A \rightarrow \B$ of the projection, we find from \eqref{eq:boundary-map} that
\begin{equation}
\delta^{r|1|s}(a_1,\dots,a_r,\boldb,a_1',\dots,a_s') = (-1)^{\circ+1} (id - \sigma \pi)
 \mu_\B^{r|1|s}(a_1,\dots,a_r,\sigma(\boldb),a_1',\dots,a_s').
\end{equation}
We usually prefer to work with $\delta[-1]: (\B/\A)[-1] \rightarrow \A$, which has degree zero and where the sign reduces to a single $-1$.

\begin{example}
Let $\A$ be an $A_\infty$-algebra, and $\PP$ an $\A$-bimodule.
One can then introduce the {\em trivial extension algebra} $\B = \A \oplus \PP$, whose structure maps are
\begin{equation} \label{eq:trivial-extension}
 \mu_\B^n =
 \begin{cases} \mu_\A^n & \text{when all entries lie in $\A$,} \\
 (-1)^{\circ+1} \mu_\PP^{r|1|n-1-r} & \text{when one entry $\boldp = b_{r+1}$ lies in $\PP$, and the others in $\A$,} \\
 0 & \text{whenever at least two entries lie in $\PP$.}
 \end{cases}
\end{equation}
The choice of signs is such that the quotient bimodule $\B/\A$ becomes equal to $\PP$. With this in mind, the obvious splitting $\sigma$ is a bimodule homomorphism, and if one uses that, then $\delta$ will vanish identically. Hence, the homotopy class $[\delta]$ is zero in the trivial extension case, for any choice of splitting.
\end{example}

\begin{lemma} \label{th:left-right}
The following diagram is homotopy commutative:
\begin{equation} \label{eq:left-right}
\xymatrix{
 & (\B/\A)[-1] \otimes_\A (\B/\A)[-1]
 \ar[dl]_{\1 \otimes \delta[-1]} \ar[dr]^{\delta[-1] \otimes \1} & \\
 (\B/\A)[-1] \otimes_\A \A \ar[dr]_{\epsilon_{right}} && \A \otimes_\A (\B/\A)[-1] \ar[dl]^{\;\;\epsilon_{left}} \\
 & (\B/\A)[-1]  &
}
\end{equation}
\end{lemma}

\proof Take the pre-homomorphism $\kappa: (\B/\A)[-1] \otimes_\A (\B/\A)[-1] \rightarrow (\B/\A)[-1]$ of degree $-1$, whose components are the maps
\begin{equation} \label{eq:kappa-homotopy}
\begin{aligned}
 & \kappa^{r|1|s}(a_1,\dots,a_r,\boldb \otimes a'_1 \otimes \cdots \otimes a'_t \otimes \boldb',
 a''_1,\dots,a''_s)
 = \\ & \qquad \qquad =
 \pi\mu_\B^{r+s+t+2}(a_1,\dots,a_r,\sigma(\boldb),a'_1,\dots,a'_t,\sigma(\boldb'),a''_1,\dots,a''_s).
\end{aligned}
\end{equation}
Then $\partial\kappa$ is precisely the difference between the two sides in \eqref{eq:left-right}. \qed

\begin{remark}
It is instructive to review how much of the given data we have used so far. The bimodule structure of $\B/\A$ contains part of $\mu_\B^*$, namely that where all inputs except one lie in $\A$, and where the output is projected to $\B/\A$. The expression for the boundary map $\delta$ takes the same kind of inputs, but this time projects the output to $\A$. To define $\kappa$, we wrote down expressions $\mu_\B^*$ where all inputs except two lie in $\A$, and where the output goes to $\B/\A$. It seems likely that there is an entire hierarchy of bimodule pre-homomorphisms, of which $\kappa$ is only the first member, and which eventually would involve all the structure of $\B$. We have no immediate use for this hierarchy, but it may be a worth while subject for future investigations.
\end{remark}

\subsection*{A bimodule computation}
We introduce a formal variable $t$ of degree $2$, which will be used to shift gradings. This is just a bookkeeping device for now, but it will assume a more central importance in the next section. For each $p \geq 1$ define an $\A$-bimodule $\TT^p$, whose underlying graded vector space is
\begin{equation} \label{eq:tt-p}
\TT^p = \bigoplus_{i_1 + \cdots + i_l = p} t^{i_1}\B[1] \otimes T(\barA[1]) \otimes t^{i_2}\B[1] \otimes \cdots \otimes T(\barA[1]) \otimes t^{i_l}\B[1],
\end{equation}
where $i_1,\dots,i_l \geq 1$. For $p = 1$, this is just $t\B[1] = \B[-1]$, and the $\A$-bimodule structure is as in \eqref{eq:b-a-sequence}. In general, we define
\begin{equation} \label{eq:tt-bimodule}
\begin{aligned}
& \mu_{\TT^p}^{0|1|0}(d_1 \otimes \cdots \otimes d_u) = \sum_{m,n} (-1)^\ast d_1 \otimes \cdots \otimes \mu_{\B}^n(d_{m+1},\dots,d_{m+n}) \otimes d_{m+n+1} \otimes \cdots \otimes d_u, \\
& \mu_{\TT^p}^{r|1|0}(a_1,\dots,a_r,d_1 \otimes \cdots \otimes d_u) = \sum_n
\mu_{\B}^{r+n}(a_1,\dots,a_r,d_1,\dots,d_n) \otimes d_{n+1} \otimes \cdots \otimes d_u,
\\
& \mu_{\TT^p}^{0|1|s}(d_1 \otimes \cdots \otimes d_u,a'_1,\dots,a'_s) = \sum_n (-1)^\ast
d_1 \otimes \cdots \otimes d_n \otimes \mu_\B^{u-n+s}(d_{n+1},\dots,d_u,a_1',\dots,a'_s),
\\
& \mu_{\TT^p}^{r|1|s}(a_1,\dots,a_r,d_1 \otimes \cdots \otimes d_u,a'_1,\dots,a'_s) =
\mu_\B^{r+u+s}(a_1,\dots,a_r,d_1,\dots,d_u,a'_1,\dots,a'_s),
\end{aligned}
\end{equation}
where $r,s>0$. Here, the $d$ variables can lie in $\barA[1]$ or in $t^i\B[1]$ for any $i>0$, and powers of $t$ are multiplied in the standard way.

\begin{lemma} \label{th:project-1}
Project $\TT^p$ to the summand where $l = p$ and $(i_1,\dots,i_l) = (1,\dots,1)$, and then project each $t\B[1]$ factor further to the quotient $t(\B/\A)[1]$, except for the rightmost one. The outcome is a quasi-isomorphism
\begin{equation} \label{eq:p-map}
\pi^p: \TT^p \longrightarrow t(\B/\A)[1]^{\otimes_\A p-1} \otimes_\A t\B[1].
\end{equation}
\end{lemma}

\proof Let's argue by induction on $p$, the case $p = 1$ being trivial. Define a finite decreasing filtration $W^\bullet \TT^p$ as follows. A tensor expression as in \eqref{eq:tt-p} lies in $W^{2i+1}\TT^p$ if the leftmost term has power $t^{i+1}$ or higher; and it lies in $W^{2i}\TT^p$ if the leftmost term either lies in $t^i\A$, or has a power strictly larger than $t^i$. By induction assumption, we then have quasi-isomorphisms
\begin{equation}
\begin{aligned}
 & W^{2i}\TT^p/W^{2i+1}\TT^p  = t^i\A[1] \otimes_\A \TT^{p-i} \\
 & \quad \xrightarrow{\htp}
 t^i\A[1] \otimes_\A t(\B/\A[1])^{\otimes_\A p-i-1} \otimes_\A t\B[1]
 && \text{for $1 \leq i \leq p-1$,} \\
 & W^{2i+1}\TT^p/W^{2i+2}\TT^p = t^{i+1}(\B/\A)[1] \otimes_\A \TT^{p-i-1}
 \\ & \quad \xrightarrow{\htp} t^{i+1}(\B/\A)[1] \otimes_\A
 t(\B/\A)[1]^{\otimes_\A p-i-2} \otimes_\A t\B[1] && \text{for $0 \leq i \leq p-2$.}
\end{aligned}
\end{equation}
Let $\CC^i$ be the cone of $\epsilon_{left}[1]: \A[1] \otimes_\A t(\B/\A[1])^{\otimes_\A p-i-1} \otimes_\A t\B[1] \rightarrow t(\B/\A[1])^{\otimes_\A p-i-1} \otimes_\A t\B[1]$, shifted up by $2i$. For each $1 \leq i \leq p-2$ there is a commutative diagram of chain complexes, with short exact columns,
\begin{equation} \label{eq:find-the-d}
\xymatrix{
 0 \ar[d] & 0 \ar[d] \\
 W^{2i+1}\TT^p/W^{2i+2}\TT^p \ar[r] \ar[d] &
 t^{i+1}(\B/\A)[1] \otimes_\A t(\B/\A)[1]^{\otimes_\A p-i-2} \otimes_\A t\B[1] \ar[d] \\
 W^{2i}\TT^p/W^{2i+2}\TT^p \ar[d] \ar[r] &
 \CC^i \ar[d]
 \\
 W^{2i}\TT^p/W^{2i+1}\TT^p \ar[r] \ar[d] &
 t^i\A[1] \otimes_\A t(\B/\A)[1]^{\otimes_\A (p-i-1)} \otimes_\A t\B[1]
 \ar[d] \\
 0 & 0
}
\end{equation}
Since the top and bottom $\rightarrow$ are quasi-isomorphisms, so is the middle one. But on the other hand, $\CC^i$ is acyclic by definition, hence so is $W^{2i}\TT^p/W^{2i+2}\TT^p$ for any $i \geq 1$. A similar but simpler reasoning applies to the bottom of the filtration. Namely, $W^{2p-2}\TT^p$ can be directly identified with the mapping cone of $\epsilon_{left}[1]: t^{p-1}\A[1] \otimes_\A t\B[1] \rightarrow t^p\B[1]$, hence is again acyclic.

Having shown that, it follows that projection $\TT^p \rightarrow \TT^p/W^2\TT^p = t(\B/\A)[1] \otimes_\A \TT^{p-1}$ is a quasi-isomorphism. Now $\pi^p$ can be thought of as the composition of that projection and $\1 \otimes \pi^{p-1}$, hence is itself a quasi-isomorphism by induction assumption. \qed

\subsection*{Curvature terms and inverse limits}
Given a left $\A$-module $\NN$ and a right $\A$-module $\MM$, we define a chain complex $T$, which is a kind of enlarged tensor product. Discussion of the significance of this is again deferred to the next section. The underlying graded vector space is
\begin{equation} \label{eq:t}
T = \prod_l \prod_{i_1,\dots,i_l} \NN \otimes T(\barA[1]) \otimes t^{i_1}\B[1] \otimes T(\barA[1])
\cdots \otimes t^{i_l}\B[1] \otimes T(\barA[1]) \otimes \MM,
\end{equation}
where $l \geq 0$ and $i_k \geq 1$; and the differential is
\begin{equation} \label{eq:partial-t}
\begin{aligned}
& \partial(\mathbf{m} \otimes d_1 \otimes \cdots \otimes d_u \otimes \mathbf{n}) = \\
& = \sum_s \mu_\NN^{1|s}(\mathbf{m},d_1,\dots,d_s) \otimes d_{s+1} \otimes \cdots \otimes d_u
\otimes \mathbf{n} \\
& + \sum_r (-1)^\ast \mathbf{m} \otimes d_1 \otimes \cdots \otimes d_{u-r} \otimes \mu_\MM^{r|1}(
d_{u-r+1},\dots,d_u,\mathbf{n}) \\
& + \sum_{m,n} (-1)^\ast \mathbf{m} \otimes \cdots \otimes d_1 \otimes \cdots \otimes d_m \otimes \mu_{\B}^n(d_{m+1},\dots,d_{m+n}) \otimes d_{m+n+1} \otimes \cdots \otimes d_u \otimes \mathbf{n} \\ & + \sum_m (-1)^\ast \mathbf{m} \otimes \cdots \otimes d_{m} \otimes te
\otimes d_{m+1} \cdots \otimes \mathbf{n}.
\end{aligned}
\end{equation}
Here, it is understood that those $d_k$ which lie in $t^i\B[1]$, $i>0$, act trivially on $\NN$ and $\MM$ in the first two lines of \eqref{eq:partial-t}. $T$ comes with a complete decreasing filtration $F^\bullet T$ by the total power of $t$. Most terms in \eqref{eq:partial-t} actually preserve that power, except for those in the last line which insert $te \in t\B[1]$, hence raise the total power by one. Obviously, passing to the graded spaces $F^pT/F^{p+1}T$ just kills those terms. By comparing this with \eqref{eq:tt-bimodule} and applying Lemma \ref{th:project-1} one sees that
\begin{equation} \label{eq:quotients}
F^pT/F^{p+1}T = \begin{cases} \NN \otimes_\A \MM & p = 0, \\
\NN \otimes_\A \TT^p \otimes_\A \MM \htp \NN \otimes_\A t(\B/\A)[1]^{\otimes_\A p-1} \otimes_\A
t\B[1] \otimes_\A \MM
& p>0. \end{cases}
\end{equation}
Consider the $\A$-module homomorphism
\begin{multline} \label{eq:putative-boundary}
 \rho:  t(\B/\A)[1]^{\otimes_\A p-1} \otimes_\A t\B[1] \otimes_\A \MM
 \xrightarrow{\text{\it projection}} t(\B/\A)[1]^{\otimes_\A p} \otimes_\A \MM \\
 \xrightarrow{\1 \otimes \xi} t(\B/\A)[1]^{\otimes_\A p} \otimes_\A t\A[1] \otimes_\A \MM
 \xrightarrow{\text{\it inclusion}} t(\B/\A)[1]^{\otimes_\A p} \otimes_\A t\B[1] \otimes_\A \MM,
\end{multline}
where $\xi$ is the map from \eqref{eq:quasi-inverse}, thought of as a degree $1$ homomorphism $\MM \rightarrow t\A[1] \otimes_\A \MM$. Let $R$ be the cone of $\1 \otimes \rho$, where the left factor is the identity of $\NN$. This fits into a diagram of chain complexes of the same form as \eqref{eq:find-the-d},
\begin{equation} \label{eq:find-the-d-2}
\xymatrix{
 0 \ar[d] & 0 \ar[d] \\
 F^{p+1}T/F^{p+2}T \ar[r] \ar[d] &
 \NN \otimes_\A t(\B/\A)[1]^{\otimes_\A p} \otimes_\A t\B[1] \otimes_\A \MM
 \ar[d] \\
 F^p T/F^{p+2} T \ar[d] \ar[r] &
 R \ar[d]
 \\
 F^pT/F^{p+1}T \ar[r] \ar[d] &
 \NN \otimes_\A t(\B/\A)[1]^{\otimes_\A p-1} \otimes_\A t\B[1] \otimes_\A \MM
 \ar[d] \\
 0 & 0
}
\end{equation}
where the top and bottom $\rightarrow$, hence also the middle one, are quasi-isomorphisms. Similarly but in a more direct way, $T/F^2T$ can be identified with the mapping cone of $\1 \otimes \rho$, where now $\rho = (\text{\it inclusion}) \circ \xi: \NN \rightarrow t\B[1] \otimes_\A \MM$.
If we consider the spectral sequence associated to $F^\bullet T$, then \eqref{eq:quotients} determines the $E_1$ page, and \eqref{eq:find-the-d-2} the differential on that page. We will use that knowledge to prove:

\begin{lemma} \label{th:finite-approx}
Project $T/F^{p+1}T$ to the summand where $l = p$ and $(i_1,\dots,i_l) = (1,\dots,1)$, and then project each $t\B[1]$ factor further to $t(\B/\A)[1]$. The outcome is a quasi-isomorphism
\begin{equation} \label{eq:proj-f}
T/F^{p+1}T \longrightarrow \NN \otimes_\A t(\B/\A)[1]^{\otimes_\A p} \otimes_\A \MM.
\end{equation}
\end{lemma}

\proof We start with a particularly simple special case, which is when $\B = \A \oplus \PP$ is a trivial extension algebra. In that case, one can write the cohomology of \eqref{eq:quotients} as
\begin{equation} \label{eq:columns-2}
 H(F^p T/F^{p+1} T) = \begin{cases} H(\NN \otimes_\A \MM) & p = 0, \\
 H(\NN \otimes_\A t\PP[1]^{\otimes_\A (p-1)} \otimes_\A \MM)[-1] \\ \qquad \oplus\;
 H(\NN \otimes_\A t\PP[1]^{\otimes_\A p} \otimes_\A \MM) \qquad & p>0.
 \end{cases}
\end{equation}
With respect to this isomorphism, the boundary operator of the left column in \eqref{eq:find-the-d-2} is the map induced by \eqref{eq:putative-boundary}. Concretely, for $p>0$ it is an isomorphism from the second summand in \eqref{eq:columns-2} into the first one, while for $p = 0$ it is the inclusion into the first summand.

Consider the spectral sequence associated to the induced filtration of $T/F^{p+1}T$, for some fixed $p$. On the $E_2$ page of that spectral sequence, only the $p$-th column is nonzero, and that is reduced to the second summand of \eqref{eq:columns-2} which is $H(\NN \otimes_\A t\PP[1]^{\otimes_\A p} \otimes_\A \MM)$, hence the same as the cohomology of the right hand side of \eqref{eq:proj-f}. An obvious comparison argument then proves the desired result, in the trivial extension case.

To do the general case, write $\PP = \B/\A$ and choose a splitting $\sigma: \PP \rightarrow \B$ of the projection. Write each $t^i\B[1]$ factor in \eqref{eq:t} accordingly as $t^i\A[1] \oplus t^i\PP[1]$, and define a finite increasing filtration $V_\bullet(T/F^{p+1}T)$ by counting the total number of $\PP$ factors. Concretely,
\begin{equation} \label{eq:v-filtration}
\begin{aligned}
& V_{-1}(T/F^{p+1}T) = 0, \\
& V_0(T/F^{p+1}T) = \bigoplus_{i_1 +\cdots +i_l \leq p}
 \NN \otimes T(\barA[1]) \otimes t^{i_1}\A[1] \otimes \cdots \otimes t^{i_l}\A[1] \otimes T(\barA[1]) \otimes \MM, \\ & \cdots \\ &
 V_p(T/F^{p+1}T)/V_{p-1}(T/F^{p+1}T) =
 \NN \otimes T(\barA[1]) \otimes t\PP[1] \otimes \cdots \otimes t\PP[1] \otimes T(\barA[1]) \otimes \MM, \\
 & V_p(T/F^{p+1}T) = T/F^{p+1}T.
\end{aligned}
\end{equation}
The map \eqref{eq:proj-f} is just projection to the last quotient in \eqref{eq:v-filtration}. Passing to the associated graded space of this filtration just means killing all terms in \eqref{eq:partial-t} which decrease the number of $\TT$ factors, which is the same as passing from a given $\B$ to the trivial extension algebra $\A \oplus\PP$. We know that in that case, \eqref{eq:proj-f} is a quasi-isomorphism, and the general result follows from that by a standard filtration argument. \qed

To be able to take the limit $p \rightarrow \infty$, we also need to know going from $p+1$ to $p$ affects the right hand side of \eqref{eq:proj-f}.

\begin{lemma} \label{th:htp-commute}
Let $\delta$ be the boundary map of \eqref{eq:b-a-sequence}, for some choice of $\sigma$.
Then the following diagram of chain maps is homotopy commutative:
\begin{equation} \label{eq:proj-inc}
\xymatrix{
 \ar[d]^-{\text{\it projection}} T/F^{p+2}T \ar[rr]^-{\text{\it projection}} && T/F^{p+1}T
 \ar[d]_-{\text{\it projection}}
 \\
 \NN \otimes_\A t(\B/\A)[1]^{\otimes_\A p+1} \otimes_\A \MM
 \ar[dr]_-{\1 \otimes \cdots \otimes \1 \otimes \delta[-1] \otimes \1\hspace{3em}} &  &
 \ar[dl]^-{\hspace{2em} \1 \otimes \cdots \otimes \1 \otimes \xi} \NN \otimes_\A t(\B/\A)[1]^{\otimes_\A p} \otimes_\A \MM
 \\
 & \hspace{-3em} \NN \otimes_\A t(\B/\A)[1]^{\otimes_\A p} \otimes_\A \A \otimes_\A \MM \hspace{-3em}&
}
\end{equation}
\end{lemma}

\proof Take $id-\sigma\pi$, think of it as a degree $-1$ pre-homomorphism $t\B[1] \rightarrow \A$, and use that to define a map $T/F^{p+2}T \rightarrow \NN \otimes_\A t(\B/\A)[1]^{\otimes_\A p} \otimes_\A t\B[1] \otimes_\A \MM \rightarrow \NN \otimes_\A t(\B/\A)[1]^{\otimes_\A p} \otimes_\A \A \otimes_\A \MM$. This provides the desired chain homotopy between the two sides of \eqref{eq:proj-inc}. \qed

The consequence is that the cohomology of $T$ fits into a short exact sequence
\begin{multline} \label{eq:h-t}
\qquad 0 \longrightarrow \underleftarrow{\mathit{lim}}^1_p\, H(\NN \otimes_\A (\B/\A)[-1]^{\otimes_\A p} \otimes_\A \MM)[-1] \longrightarrow H(T) \\ \longrightarrow
\underleftarrow{\mathit{lim}}_p\, H(\NN \otimes_\A (\B/\A)[-1]^{\otimes_\A p} \otimes_\A \MM) \longrightarrow 0 \qquad
\end{multline}
where the inverse limit, and its derived functor, are both formed with respect to $\delta$, in the same sense as in \eqref{eq:proj-inc}. Note that by a trivial application of the Mittag-Leffler condition, the derived term vanishes if each $H(\NN \otimes_\A t(\B/\A)[1]^{\otimes_\A p} \otimes_\A \MM)$ is finite-dimensional \cite[p.\ 83]{weibel}.

\section{Categories of modules\label{sec:categories}}

This section introduces the second main object, namely the curved $A_\infty$-algebra $\DD$. Having defined that, we will study the relation between the categories of modules over $\A$ and over $\DD$, leading to the main result of the paper (Theorem \ref{th:main}).

\subsection*{The boundary map as a natural transformation}
First, we want to reinterpret Lemma \ref{th:left-right} in terms of convolution functors. Let $\VV = mod(\A)$. The bimodule $(\B/\A)[-1]$ defines a dg functor $\Phi_{(\B/\A)[-1]}: \VV \rightarrow \VV$, and the homomorphism $\delta[-1]: (\B/\A)[-1] \rightarrow \A$ gives a natural transformation $\Phi_{\delta[-1]}: \Phi_{(\B/\A)[-1]} \rightarrow \Phi_{\A}$. To avoid technical problems stemming from the failure of \eqref{eq:a-fun} to be strictly commutative, we pass to the cohomology level, where the relevant structures are \begin{equation}
\begin{aligned}
& V = H^0(\VV), \\
& F = H^0(\Phi_{(\B/\A)[-1]}): V \longrightarrow V, \\
& T = [\Phi_{\delta[-1]}]: F \longrightarrow H^0(\Phi_{\A}) \iso Id.
\end{aligned}
\end{equation}
$T$ is ambidextrous. This is because the induced natural transformations $R_FT$, $L_FT$ from \eqref{eq:two-sides} are given by the two sides of the following diagram, which is homotopy commutative by Lemma \ref{th:left-right}:
\begin{equation} \label{eq:left-right-2}
\xymatrix{
 & (\B/\A)[-1] \otimes_\A (\B/\A)[-1] \otimes_\A \MM
 \ar[dl]_{\1 \otimes \delta[-1] \otimes \1\;\;} \ar[dr]^{\;\;\delta[-1] \otimes \1 \otimes \1} & \\
 (\B/\A)[-1] \otimes_\A \A \otimes_\A \MM \hspace{-4em}
 \ar[dr]_{\1 \otimes \epsilon_{left} \htp \epsilon_{right} \otimes  \1 \hspace{4em}} &&
 \hspace{-4em} \A \otimes_\A (\B/\A)[-1] \otimes_\A \MM \ar[dl]^{\;\;\epsilon_{left} \otimes \1} \\
 & (\B/\A)[-1] \otimes_\A \MM.  &
}
\end{equation}

\subsection*{A curved $A_\infty$-algebra}
Consider the graded vector space
\begin{equation} \label{eq:filtered-morphisms}
 \DD = \A \oplus t\B[[t]] \subset \B[[t]]
\end{equation}
where, as before, $t$ is a formal variable of degree two. In words, elements of $\DD^k$ are formal power series in $t$, of total degree $k$, constrained by asking that the constant ($t^0$) term should lie in $\A$. The $A_\infty$-structure maps $\mu_\DD^s$, $s>0$, are the $t$-linear extensions of those on $\B$. Additionally, we introduce a curvature term
\begin{equation} \label{eq:curvature-term}
 \mu_\DD^0 = te \in \DD^2,
\end{equation}
where $e$ is the identity on $\A$. This satisfies the extended $A_\infty$-equations, which are the same as \eqref{eq:a} but allowing terms with $n = 0$ as well.

$\DD$ comes with its $t$-adic filtration $F^\bullet\DD$, given by $F^0\DD = \DD$, $F^p = t^p\B[[t]]$ for $p>0$. This is a complete decreasing filtration, compatible with the $A_\infty$-structure, and the curvature term \eqref{eq:curvature-term} is small, meaning that it lies in $F^1\DD$. These are precisely the conditions required to make $\DD$ into a {\em filtered curved $A_\infty$-algebra}, which is a special case of the algebraic structures studied in \cite{fooo}. Note that in spite of the obvious $t$-linearity of composition maps, $\DD$ is still considered to be defined over $\K$, and not $\K[[t]]$. However, when making any constructions involving it, the $t$-adic topology needs to be taken into account.

\subsection*{A category of modules}
The naive definition of $\DD$-module would be a graded vector space $\MM$, together with a structure map $\mu_\MM: T(\DD[1]) \otimes \MM \rightarrow \MM[1]$, which satisfies an appropriately extended version of the $A_\infty$-module equation. Such modules no longer have cohomology, since the first of the structure equations is
\begin{equation}
\mu^1_\MM(\mu^1_\MM(\mathbf{m})) + \mu^2_\MM(\mu^0_\DD,\mathbf{m}) +
(-1)^{|m|} \mu^2_\MM(\mathbf{m},\mu^0_\DD) = 0.
\end{equation}
We will need to modify the general definition somewhat, in order to adapt it to our particular context. A {\em torsion $\DD$-module} is one which admits a finite decreasing filtration $F^\bullet\MM$, such that:
\begin{equation}
\parbox{36em}{
 If $d_1 \in t^{i_1}\DD, \dots, d_r \in t^{i_r}\DD$ and $\mathbf{m} \in F^j\MM$, then
 $\mu_\MM^{r|1}(d_1,\dots,d_r,\mathbf{m}) \in F^{i_1 + \cdots + i_r +j}\MM$.
}
\end{equation}
Projection $\DD \rightarrow \DD/t\B[[t]] = \A$ allows us to pull back $\A$-modules to $\DD$. Obviously, a $\DD$-module is such a pullback if and only if $t\B[[t]] \subset \DD$ acts trivially on it. In these terms, a torsion $\DD$-module is one which has a filtration such that the graded pieces $F^p\MM/F^{p+1}\MM$ are pulled back from $\A$.

In line with the general philosophy concerning $\DD$, we only want to allow module pre-homomorphisms which are continuous in the $t$-adic topology. In the case of torsion modules, where the topology is discrete, this means that a pre-homomorphism of degree $k$ is given by a map $\phi: T(\bar{\DD}[1]) \otimes \MM \rightarrow \NN[k]$ with the following property:
\begin{equation} \label{eq:continuous}
\parbox{36em}{
There is some $q \gg 0$ such that if $d_1 \in t^{i_1}\DD, \dots, d_r \in t^{i_r}\DD$ with $d_1 + \cdots + d_r > q$, then $\phi^{r|1}(d_1,\dots,d_r,\mathbf{m}) = 0$.}
\end{equation}
The differential and composition are defined in the standard way, taking $\mu^0_\DD$ into account. We denote the resulting dg category of torsion modules by $\WW = modt(\DD)$, and the pullback dg functor by $\GG: \VV \rightarrow \WW$. On the cohomology level, we have
\begin{equation}
\begin{aligned}
& W = H^0(\WW), \\
& G = H^0(\GG): V \longrightarrow W.
\end{aligned}
\end{equation}
As a consequence of \eqref{eq:continuous}, the cone of any degree $1$ cocycle in $\WW$ is again a torsion module, which means that $\WW$ is pre-triangulated. As a consequence, $W$ is triangulated, and $G$ is an exact functor. Moreover, because of the existence of filtrations, the objects in the image of $G$ generate $W$.

A parallel discussion applies to right modules. In fact, the chain complex $T$ discussed in Section \ref{sec:subalgebra} is just the tensor product $\MM \otimes_\DD \NN$, where $\MM$ is a right $\A$-module and $\NN$ a left $\A$-module, both being pulled back to $\DD$ (and where the tensor product is really a topological, which means $t$-adically completed, one). Therefore, \eqref{eq:h-t} is a statement about the behaviour of modules under pullback. Rather than elaborate on that, we will consider the dual construction involving $hom$s instead of tensor products.

\subsection*{Morphisms of pullback modules}
Let $\MM$ and $\NN$ be two $\A$-modules. Pull both back to $\DD$, and consider the chain complex $C = hom_\WW(\MM,\NN)$. Explicitly,
\begin{equation}
 C = \bigoplus_l \bigoplus_{i_1, \dots, i_l} hom(T(\barA[1]) \otimes t^{i_1}\B[1] \otimes T(\barA[1]) \otimes \cdots \otimes t^{i_l}\B[1] \otimes T(\barA[1]) \otimes \MM,\NN)
\end{equation}
where the $hom$ on the right hand side is just the graded vector space of linear maps. The fact that we have a direct sum instead of a product comes from \eqref{eq:continuous}. The differential is
\begin{equation} \label{eq:d-c}
\begin{aligned}
& (\partial\phi)^{r|1}(d_1,\dots,d_r,\mathbf{m}) =
 \\ &
 \sum_{m}
 (-1)^{|\phi| \ast} \
 \mu_\NN^{m|1}(d_1,\dots,d_m,\phi^{r-m|1}(d_{m+1},\dots,d_r,\mathbf{m}))
 \\ + &
 \sum_{m}
 (-1)^{|\phi|+1+\ast}
 \phi^{m|1}(d_1,\dots,d_m,\mu_\MM^{r-m|1}(d_{m+1},\dots,d_r,\mathbf{m}))
 \\ + &
 \sum_{m,n}
 (-1)^{|\phi|+1+\ast}
 \phi^{r-n+1|1}(d_1,\dots,d_m,\mu_\B^n(d_{m+1},\dots,d_{m+n}),d_{m+n+1},\dots,\mathbf{m})
 \\ + &
 \sum_m
 (-1)^{|\phi|+1+\ast}
 \phi^{r+1|1}(d_1,\dots,d_m,te,d_{m+1},\dots,d_r,\mathbf{m}),
\end{aligned}
\end{equation}
where it is understood that $t\B[[t]] \subset \DD$ acts trivially on both $\MM$ and $\NN$. This means that the first line of \eqref{eq:d-c} vanishes unless $d_1,\dots,d_m \in \bar\A[1]$, and that the second line vanishes unless $d_{m+1},\dots,d_r \in \bar\A[1]$.
Let $F_\bullet C$ be the increasing $t$-adic filtration of $C$; $F_p C$ consists of those maps which vanish if the total power of $t$ involved in the argument is $>p$. Passing to the graded spaces $F_pC/F_{p-1}C$ just kills the curvature term, which is the last line in \eqref{eq:d-c}. Comparing this to \eqref{eq:tt-bimodule} shows that
\begin{equation}
F_pC/F_{p-1}C = \begin{cases} hom_\VV(\MM,\NN) & p = 0, \\
hom_\VV(\TT^p \otimes_\A \MM,\NN) \htp
hom_\VV(t(\B/\A)[1]^{\otimes_\A p-1} \otimes_\A t\B[1] \otimes_\A \MM,\NN) & p>0,
\end{cases}
\end{equation}
which is the analogue of \eqref{eq:quotients} in our context. Proceeding as in \eqref{eq:find-the-d-2}, one can show that $F_{p+1}C/F_{p-1}C$ is the mapping cone of the chain map induced by the module homomorphism $\rho$ from \eqref{eq:putative-boundary}, or its simplified version for $p = 1$. From there, following exactly the same path as in Lemma \ref{th:finite-approx}, one finds that the inclusions
\begin{equation} \label{eq:include-p}
\iota_p: hom_\VV(t(\B/\A[1])^{\otimes_\A p} \otimes_\A \MM,\NN) \longrightarrow F_pC = F_p hom_\WW(\MM,\NN)
\end{equation}
are quasi-isomorphisms. On the cohomological level, this means that
\begin{equation} \label{eq:coho-direct}
\underrightarrow{lim}_p\, Hom_V((\B/\A)[-1]^{\otimes_\A p} \otimes_\A \MM,\NN) \iso
Hom_W(\MM,\NN).
\end{equation}
In principle, one could follow the same strategy as in Lemma \ref{th:htp-commute} to prove that the maps which occur in this direct system are induced by $\delta[-1]$. However, it turns out that there is an alternative approach, which leads to the same conclusion in a form which is better adapted for our uses.

\subsection*{The multiplicative structure}
Define $\psi \in hom_\WW(\MM,(\B/\A)[-1] \otimes_\A \MM)$ as follows: for $b \in \B$ and $a_2,\dots,a_r \in \barA$, $\psi^{r|1}(tb,a_2,\dots,a_r,\mathbf{m}) = \pi(b) \otimes a_2 \otimes \cdots \otimes a_r \otimes \mathbf{m}$; and all other components are zero. It is easy to see that $\partial\psi = 0$. Composition with $\psi$ induces a map $hom_\WW((\B/\A)[-1] \otimes_\A \MM,\NN) \rightarrow hom_\WW(\MM,\NN)$ which raises the $t$-adic filtration level by one, and which fits into a commutative diagram
\begin{equation} \label{eq:multiply-psi}
\xymatrix{hom_\VV(t(\B/\A)[1]^{\otimes_\A p} \otimes_\A (\B/\A)[-1] \otimes_\A \MM,\NN)
\ar@{=}[r] \ar[d]^{\htp}_{\iota_p} & \ar[d]^{\htp}_{\iota_{p+1}}
hom_\VV(t(\B/\A)[1]^{\otimes_\A p+1} \otimes_\A \MM,\NN) \\
F_p hom_\WW((\B/\A)[-1] \otimes_\A \MM,\NN) \ar[r]^-{\psi} &
F_{p+1} hom_\WW(\MM,\NN)
}
\end{equation}
Passing to the direct limit in $p$, we see that composition with $\psi$ is a quasi-isomorphism. By standard categorical nonsense, $[\psi]$ itself is an isomorphism in the cohomological category $W$. Moreover, the combination of various diagrams of type \eqref{eq:multiply-psi} yields a commutative diagram
\begin{equation} \label{eq:multiply-psi-2}
\xymatrix{ & \hspace{-5em} hom_\VV(t(\B/\A)[1]^{\otimes_\A p} \otimes_\A \MM,\NN) \ar[dl]^{=}_{\GG}
\ar[dr]^{\iota_p}_{\htp}
\hspace{-5em} & \\
F_0 hom_\WW((\B/\A)[-1]^{\otimes_\A p} \otimes_\A \MM,\NN) \ar[rr] &&
F_p hom_\WW(\MM,\NN)
}
\end{equation}
where the bottom $\rightarrow$ is successive composition with the homomorphisms $\xi$ associated to the modules $\MM, (\B/\A)[-1] \otimes_\A \MM, \dots, (\B/\A)[-1]^{\otimes_\A p-1} \otimes_\A \MM$.

Additionally, we have a pre-homomorphism $\nu \in hom_\WW(\MM,\A \otimes_\A \MM)$ of degree $-1$, given by a similar formula as $\psi$, but using the splitting $\sigma$ and the corresponding projection $id-\sigma\pi: \B \rightarrow \A$. This means that for $b \in \B$ and $a_2,\dots,a_r \in \barA$, $\nu^{r|1}(tb,a_2,\dots,a_r,\mathbf{m}) = (id-\sigma\pi)(b) \otimes a_2 \otimes \cdots \otimes a_r \otimes \mathbf{m}$. This has the property that
\begin{equation}
\partial \nu = \xi - (\delta[-1] \otimes \1) \psi
\end{equation}
where both $\xi$ and $\delta[-1] \otimes \1$ come from $\VV$ through the pullback functor. It follows that $\xi$ is the inverse up to homotopy of $\GG(\epsilon (\delta[-1] \otimes \1)) \in hom_\WW((\B/\A)[-1] \otimes_\A \MM,\MM)$. On the cohomological level, this means first of all that the elements defining the natural transformation $T$ become isomorphisms under $G: V \rightarrow W$. Moreover, thanks to \eqref{eq:multiply-psi-2}, the isomorphism in \eqref{eq:coho-direct} can be defined by applying $G$ and then multiplying with the inverses of $T$, as required in Lemma \ref{th:dg-quotient}. The other assumptions in that Lemma having been checked earlier, we get:

\begin{theorem} \label{th:main}
$\WW = modt(\DD)$ is quasi-equivalent to the dg quotient of $\VV = mod(\A)$ by the full dg subcategory $\VV_{nil}$ associated to the natural transformation $T$. In particular, the triangulated category $W = H^0(modt(\DD))$ is the localization of $V = H^0(mod(\A))$ along $T$. \qed
\end{theorem}

\subsection*{Complements of divisors revisited}
Let's return to one of the two examples from Section \ref{sec:loc}, although in somewhat reduced generality. Take a smooth affine variety $X$ over $\K = \C$, and $Y \subset X$ a hypersurface, whose defining function is $s \in \C[X]$. We set $\A = \C[X]$, and define a commutative dga $\B = \A[\epsilon]$, where $|\epsilon|=-1$ and $\partial\epsilon = s$, giving $H(\B) = \C[Y]$ (the advantage of using this resolution for $Y$ is that restriction of functions turns into the inclusion $\A \rightarrow \B$, which fits well into our general framework). Turn $\A$ and $\B$ into $A_\infty$-algebras in the obvious way. We have $\B/\A = \A[1]$, so the boundary homomorphism of \eqref{eq:b-a-sequence} is a bimodule homomorphism $\A \rightarrow \A$ of degree zero, which is just given by $s$ itself. Let $V^{perf} \subset V$ be the smallest thick triangulated subcategory containing the free module $\A$, and $W^{perf} \subset W$ its image under $G: V \rightarrow W$. The first of these is just $V^{perf} = D^b(X)$, and by Theorem \ref{th:main} we find that $W^{perf} = D^b(U)$, where $U = X \setminus Y$.

On the other hand, one can write $\DD = \A \oplus t\B[[t]] = \A[[t]][\gamma]$, where $\gamma = \epsilon t$ has degree one (as usual, $\mu^0_\DD = t$, and the other nontrivial terms $\mu^1, \mu^2$ are inherited from the dg structure of $\B$). Let's apply {\em Koszul duality relative to $\A$}. Namely, take the $\DD$-module $\A$, and consider its endomorphism dga, which turns out to be quasi-isomorphic to the commutative dga $\CC = \A[\theta][g]$, where $|\theta| = -1$, $|g| = 0$, and the differential is $\partial\theta = 1-gs$. In particular, $H(\CC) = \A/(gs-1) = \C[U]$. Convolution provides an embedding of $W^{perf} \subset W$ into $H^0(mod(\CC)) \iso H^0(mod(\C[U]))$. In this way, the relationship with $U$ emerges in a slightly more transparent way. It may be possible to use a suitably sheafified version of this to address the general quasi-projective case.

\section{Hochschild homology}
This section contains the closed string analogues of the computations from Section \ref{sec:subalgebra}. The ``closed string'' terminology is loosely inspired by the role of Hochschild homology in string theory, but one can also take it in quite a naive sense, meaning that we deal with tensor products in the form of closed chains. The target of our discussion is Theorem \ref{th:hochschild}, which is the counterpart of \eqref{eq:h-t}.

\subsection*{Cyclic tensor products}
Fix $A_\infty$-algebras $\A_1,\dots,\A_l = \A_0$. Suppose moreover that for each $1 \leq i \leq l$ we have a bimodule $\PP_i$ over $(\A_{i-i},\A_i)$. The {\em cyclic tensor product} of these bimodules, denoted by
$
  Z= \PP_1 \otimes_{\A_1} \PP_2 \otimes_{\A_2} \cdots \otimes_{\A_{l-1}} \PP_l \otimes_{\A_l} \cycl
$,
is the following chain complex:
\begin{equation} \label{eq:cyclic-tensor}
\begin{aligned}
 & Z  = \PP_1 \otimes T(\barA_1[1]) \otimes \PP_2 \otimes \cdots \otimes T(\barA_{m-1}[1]) \otimes \PP_l \otimes T(\barA_m[1]), \\[0.5em]
 & \partial(\boldp_1 \otimes a_{1,1}
 \otimes \cdots \otimes a_{1,u_1} \otimes \boldp_2 \otimes \cdots \otimes
 \boldp_l \otimes a_{l,1} \otimes \cdots \otimes a_{l,u_l}) \\ & =
 \sum_{i,m,n} (-1)^\ast \, \boldp_1 \otimes \cdots \otimes \boldp_i \otimes \cdots a_{i,m} \otimes
 \mu_{\A_i}^n(a_{i,m+1},\dots,a_{i,m+n}) \otimes \cdots \boldp_{i+1} \otimes \cdots \otimes a_{l,u_l} \\
 & +
 \sum_{i,m,n} (-1)^\ast \, \boldp_1 \otimes \cdots a_{i,m} \otimes
 \mu_{\PP_i}^{u_i-m|1|n}(a_{i,m+1},\dots,\boldp_i,\dots,a_{i+1,n}) \otimes
\cdots
 \otimes a_{l,u_l} \\
 & +
 \sum_{m,n} (-1)^{\#}
 \mu_{\PP_1}^{u_l-m|1|n}(a_{l,m+1},\dots,a_{l,u_l},\boldp_1,a_{1,1},\dots
 a_{1,n}) \otimes \cdots \otimes a_{l,m},
\end{aligned}
\end{equation}
where $\# = (\|a_{l,m+1}\| + \cdots + \|a_{l,u_l}\|)(|\boldp_1| + \cdots + \|a_{l,m}\|)$. The last term in the differential explains the name ``cyclic''. Given bimodule pre-homomorphisms $\phi_i: \PP_i \rightarrow \QQ_i$ for $1 \leq i \leq l$, one defines a map $f = \phi_1 \otimes \cdots \otimes \phi_l \otimes \cycl$ between the respective cyclic tensor products, as follows:
\begin{equation} \label{eq:cycl-map}
\begin{aligned}
 & f(\boldp_1 \otimes a_{1,1}
 \otimes \cdots \otimes a_{1,u_1} \otimes \boldp_2 \otimes \cdots \otimes
 \boldp_l \otimes a_{l,1} \otimes \cdots \otimes a_{l,u_l}) \\
 & \qquad = (-1)^{\#+\S} \phi_1^{u_l-m_l|1|n_1}(a_{l,m_l+1},\dots,a_{l,u_l},\boldp_1,a_{1,1},\dots,
 a_{1,n_1}) \otimes a_{1,n_1+1} \otimes \cdots \otimes a_{1,m_1} \\ & \qquad \otimes
 \phi_2^{u_1-m_1|1|n_2}(a_{1,m_1+1},\dots,a_{1,u_1},\boldp_2,a_{2,1},\dots,
 a_{2,n_2}) \otimes a_{2,n_2+1} \otimes \cdots \otimes a_{2,m_2} \\
 & \qquad \dots \\
 & \qquad \otimes \phi_l^{u_{l-1}-m_{l-1}|1|n_l}(a_{l-1,m_{l-1}+1}, \dots,
 a_{l-1,u_{l-1}}, \boldp_l, a_{l,1}, \dots, a_{l,n_l}) \otimes a_{l,n_l+1} \otimes \cdots
 \otimes a_{l,m_l}.
\end{aligned}
\end{equation}
Here $\# = (\|a_{l,m_l+1}\| + \cdots + \|a_{l,u_l}\|)(|\boldp_1| + \cdots + \|a_{l,m_l}\|)$ is the same permutation Koszul sign as before, and $\S = |\phi_2|(\|a_{l,m_l+1}\| + \cdots + |\boldp_1| + \cdots + \|a_{1,m_1}\|) + |\phi_3|(\|a_{l,m_l+1}\| + \cdots + |\boldp_1| + \cdots + |\boldp_2| + \cdots + \|a_{2,m_2}\|) + \cdots$ is the Koszul sign which arises when moving the $\phi_k$ from the left to their appropriate positions. \eqref{eq:cycl-map} is compatible with differentials, meaning in particular that if the $\phi_i$ are homomorphisms, then $f$ is a chain map. Moreover, if the $\phi_i$ are quasi-isomorphisms, then so is $f$. Finally, just like the ordinary tensor product, this construction behaves well with respect to composition of homomorphisms.

Even though our terminology is non-standard, the notion itself is not new by any means. In the simplest case, where one has only one algebra $\A = \A_0 = \A_1$ and bimodule $\PP = \PP_1$, $Z$ is the reduced Hochschild complex of $\A$ with coefficients in $\PP$, and its cohomology is called the Hochschild homology $HH(\A,\PP)$ with coefficients in $\PP$. Specializing even further to $\PP = \A$, one writes $HH(\A) = HH(\A,\A)$. More general cyclic tensor products can also be written as Hochschild homology groups, simply by noticing that
\begin{equation}
 \PP_1 \otimes_{\A_1} \PP_2 \otimes_{\A_2} \cdots \otimes_{\A_{l-1}} \PP_l \otimes_{\A_l} \cycl
 = (\PP_1 \otimes_{\A_1} \cdots \otimes_{\A_{l-1}} \PP_l) \otimes_{\A_l} \cycl.
\end{equation}

There is one case where the tensor product reduces to an ordinary non-cyclic one, namely when $\PP_1 = \MM \otimes \NN$ is the tensor product (over $\K$) of a left $\A$-module $\MM$ and right $\A$-module $\NN$. By moving the first factor to the right (with suitable Koszul signs) one then obtains an isomorphism $Z \iso \NN \otimes_{\A_1} \PP_2 \cdots \otimes_{\A_{l-1}} \PP_l \otimes_{\A_l} \MM$. In order to extend the domain of applicability of this trick, it is useful to have some form of free resolutions of bimodules. This is well-known for modules over dg algebras \cite[Section 10.12.2.4]{bl}, and the proof generalizes in a fairly direct way, but it may still be worth while to reproduce the argument:

\begin{lemma} \label{th:free}
Given any $(\A,\A')$-bimodule $\PP$, one can find a quasi-isomorphism $\QQ \rightarrow \PP$, where $\QQ$ has the following property: it admits a bounded below increasing filtration $L_\bullet \QQ$, such that each graded space $L_i\QQ/L_{i-1}\QQ$ is a direct sum of shifted copies of the bimodule $\A \otimes \A'$.
\end{lemma}

\proof Construct an increasing sequence $\PP_0 \subset \PP_1 \subset \PP_2 \cdots$ as follows. The starting point is $\PP_0 = \PP$. Homotopy classes of bimodule homomorphisms $\A \otimes \A' \rightarrow \PP_{i-1}$ are in bijective correspondence with elements of $H(\PP_{i-1})$; this takes $[\phi]$ to $[\phi^{0|1|0}(e \otimes e')]$, and its bijectivity can be reduced to the analogous fact in classical module theory by a spectral sequence argument. A parallel statement holds for direct sums of shifted copies of $\A \otimes \A'$. In particular, up to homotopy there is a unique bimodule homomorphism $\A \otimes H(\PP_{i-1}) \otimes \A' \rightarrow \PP_{i-1}$ corresponding to $Id: H(\PP_{i-1}) \rightarrow H(\PP_{i-1})$. By construction, this induces a surjective map on cohomology. Take $\PP_i$ to be the cone of that homomorphism. The union $\PP_\infty = \bigcup_i \PP_i$ is acyclic, hence if we take $\QQ = (\PP_\infty/\PP_0)[-1]$, the boundary map $\QQ \rightarrow \PP_0 = \PP$ is a quasi-isomorphism. The desired filtration is $L_i \QQ = (\PP_i/\PP_0)[-1]$. \qed

\subsection*{The Hochschild homology of $\DD$, with torsion coefficients}
In the same way as in Section \ref{sec:categories}, one can construct a dg category $bimodt(\DD)$ of torsion $\DD$-bimodules. If $\PP_1, \dots, \PP_l$ are such bimodules, we define $\PP_1 \otimes_{\DD} \PP_2 \cdots \otimes_{\DD} \PP_l \otimes \cycl$ by an appropriate generalization of \eqref{eq:cyclic-tensor}, keeping in mind that the tensor algebra $T(\bar{\DD}[1])$ should be used in its $t$-adically completed form, and that the differential has additional $\mu^0_\DD$ terms.

At the moment, the only relevant case will be when there is a single bimodule, which moreover is pulled back from $\A$. Fix such a $\PP$, write $X = \PP \otimes_{\DD} \cycl$ for the Hochschild complex formed over $\DD$, and $HH(\DD,\PP) = H(X)$ for its cohomology. $X$ carries a complete decreasing filtration $F^\bullet X$ by total powers of $t$, and the finite quotients of that filtration come with natural projections
\begin{equation} \label{eq:project-3}
 X/F^{p+1}X \longrightarrow \PP \otimes_\A t(\B/\A)[1]^{\otimes_\A p} \otimes_\A \cycl.
\end{equation}
Our first claim is that these are quasi-isomorphisms. If $\PP = \MM \otimes \NN$, shifting the $\MM$ factor to the right identifies $X$ with the complex $T$ from \eqref{eq:t}, and \eqref{eq:project-3} with \eqref{eq:proj-f}, which we know to be a quasi-isomorphism. To derive the general result from this, we note that an $\A$-bimodule quasi-isomorphism $\QQ \rightarrow \PP$ induces quasi-isomorphisms on both sides of \eqref{eq:project-3}. One can therefore replace any given $\PP$ by one of the kind considered in Lemma \ref{th:free}. Then, an obvious filtration argument applies, reducing the desired statement to the case of the bimodule $\A \otimes \A$, which is part of the previously considered situation.

The next claim is that the maps \eqref{eq:project-3} fit into a homotopy commutative diagram
\begin{equation}
\xymatrix{
 \ar[d]^-{\text{\it projection}} X/F^{p+2}X \ar[rr]^-{\text{\it projection}} && X/F^{p+1}X
 \ar[d]_-{\text{\it projection}}
 \\
 \PP \otimes_\A t(\B/\A)[1]^{\otimes_\A p+1} \otimes_\A \cycl
 \ar[dr]_-{\1 \otimes \cdots \otimes \1 \otimes \delta[-1] \otimes \cycl\hspace{3em}} &  &
\ar[dl]^-{\hspace{2em} \1 \otimes \cdots \otimes \1 \otimes \xi \otimes \cycl}
\PP \otimes_\A t(\B/\A)[1]^{\otimes_\A p} \otimes_\A \cycl
 \\
 & \hspace{-4em} \PP \otimes_\A t(\B/\A)[1]^{\otimes_\A p} \otimes_\A \A \otimes_\A \cycl \hspace{-4em}&
}
\end{equation}
Note that in both diagonal arrows, the leftmost homomorphism is the identity, which means that the entries in \eqref{eq:cycl-map} do not actually get permuted. The statement is obviously parallel to
\eqref{eq:proj-inc}, and the required chain homotopy can be constructed in the same way as in Lemma \ref{th:htp-commute}. The upshot is:

\begin{lemma} \label{th:twisted-hh}
For any $\A$-bimodule $\PP$, there is a short exact sequence
\begin{multline}
\qquad 0 \longrightarrow \underleftarrow{\mathit{lim}}^1_p\, HH(\A,\PP \otimes_\A (\B/\A)[-1]^{\otimes_\A p}) \\ \longrightarrow HH(\DD,\PP) \longrightarrow
\underleftarrow{\mathit{lim}}_p\, HH(\A,\PP \otimes_\A (\B/\A)[-1]^{\otimes_\A p})
\longrightarrow 0, \qquad
\end{multline}
where the connecting maps in the inverse limit are given by applying the homomorphism $\delta$ to one of the $\B/\A$ factors. \qed
\end{lemma}

\begin{lemma} \label{th:b-trivial}
Consider $\B$ as a torsion $\DD$-bimodule by pullback, $\DD \rightarrow \A \rightarrow \B$. Then $HH(\DD,\B) = 0$.
\end{lemma}

\proof In $bimodt(\DD)$ consider the degree $-1$ pre-homomorphism $\alpha$ from $\B$ to itself, whose only nonvanishing components are
\begin{multline}
\qquad \alpha^{r|1|s}(a_1,\dots,a_r,\mathbf{b},a'_1,\dots,a_{i-1},tb',a'_{i+1},\dots,a'_s) \\  =
- \mu_\B^{r+1+s}(a_1,\dots,a_r,\mathbf{b},a'_1,\dots,a_{i-1},b',a'_{i+1},\dots,a'_s) \qquad
\end{multline}
for $a_k, a'_j \in \bar{\A}$ and $\mathbf{b}, b' \in \B$. One computes easily that $\partial\alpha = \1$ is the identity morphism. Hence, $\B$ is isomorphic to the zero object in $H^0(bimodt(\DD))$. Hochschild homology being functorial in that category, the desired result follows immediately. \qed

\subsection*{The Hochschild homology of $\DD$, with diagonal coefficients}
Strictly speaking, the ordinary Hochschild homology $HH(\DD) = HH(\DD,\DD)$ does not fall under the definition given above, since the diagonal is not a torsion bimodule. The only necessary modification is a further $t$-adic completion of the Hochschild complex, but it still may be best to write down the resulting complex, denoted by $D$, explicitly:
\begin{equation} \label{eq:hh-dd}
\begin{aligned}
D & = \prod \A \otimes T(\barA[1]) \otimes t^{i_2}\B[1] \otimes T(\barA[1]) \otimes \cdots \otimes T(\barA[1]) \otimes t^{i_l}\B[1] \otimes T(\barA[1]) \\
 & \oplus \prod t^{i_1}\B \otimes T(\barA[1]) \otimes t^{i_2}\B[1] \otimes T(\barA[1]) \otimes \cdots \otimes T(\barA[1]) \otimes t^{i_l}\B[1] \otimes T(\barA[1]),
\end{aligned}
\end{equation}
where the product is over $l \geq 1$ and $i_1,\dots,i_l > 0$ in the second case; and over $l \geq 1$, $i_2,\dots,i_l > 0$ in the first case, containing in particular the Hochschild complex of $\A$ as the case $l = 1$. Let $G^\bullet D$ be the complete decreasing filtration of $D$ whose terms are the second summand in \eqref{eq:hh-dd}, followed by its subspace where $i_1 \geq p$ for some number $q$. Clearly, $D/G^1D$ is the Hochschild complex of $\DD$ with coefficients in $\A$, while for any $q>0$, $G^qD/G^{q+1}D$ is the analogous complex with coefficients in $\B$ (up to an even shift). The cohomology of these complexes is determined by Lemma \ref{th:twisted-hh} and Lemma \ref{th:b-trivial}, respectively. Hence:

\begin{theorem} \label{th:hochschild}
The Hochschild homology of $\DD$ fits into a short exact sequence
\begin{multline}
\qquad 0 \longrightarrow \underleftarrow{\mathit{lim}}^1_p\, HH(\A, (\B/\A)[-1]^{\otimes_\A p}) \\ \longrightarrow HH(\DD) \longrightarrow
\underleftarrow{\mathit{lim}}_p\, HH(\A,(\B/\A)[-1]^{\otimes_\A p})
\longrightarrow 0, \qquad
\end{multline}
where the connecting maps are as in Lemma \ref{th:twisted-hh}. \qed
\end{theorem}

\section{Symplectic geometry}
This section is a concise discussion of the original geometric motivation for our algebraic constructions, mostly following \cite{seidel06}. We also consider one nontrivial example, which is based on the known validity of homological mirror symmetry for $\CP{2}$ and for elliptic curves (as in \cite{polishchuk-zaslow98,polishchuk98,auroux-katzarkov-orlov04}, even though our presentation is actually a scaled-down version of \cite{seidel03b}).

\subsection*{Semisimple base rings}
Instead of a field $\K$, we will need to work over a semisimple ring of the form $R = \K^d = \K e_1 \oplus \cdots \oplus \K e_d$, where $e_i^2 = e_i$, $e_ie_j = 0$ for $i \neq j$. For instance, an $A_\infty$-algebra $\A$ over $R$ is the same as an $A_\infty$-category with $d$ ordered objects. Basis elements of $T(\bar{\A}[1])$, where the tensor product is now taken over $R$, can be thought of as composable chains of morphisms in that category. Correspondingly, when taking cyclic tensor products, one should use closed composable chains. To be precise, write $M_{\text{\it diag}} = \bigoplus_i e_i M e_i$ for the diagonal part of any $R$-bimodule $M$. Then the Hochschild complex of $\A$ is given by $(\A \otimes_R T(\bar{\A}[1]))_{\text{\it diag}}$, and similarly in the more general context of \eqref{eq:cyclic-tensor}. With this in mind, the main results, and their proofs, remain as before.

\subsection*{A conjectural dictionary}
Let $\pi: E \rightarrow D$ be an exact symplectic Lefschetz fibration over a disc, together with a trivialization of the canonical bundle of $E$. For the precise definitions, see for instance \cite[Chapter 3]{seidel04}. Fix some $z \in \partial D$, and let $M = E_z$ be the fibre at that point. $M$ is an exact symplectic manifold with contact type boundary, and comes with an induced trivialization of its canonical bundle, hence has a well-defined Fukaya category $\F(M)$, which is an $A_\infty$-category over $\K$. We will assume that $\F(M)$ is strictly unital, which is not true with the most common definition, but can always be achieved by passing from the given category to a quasi-isomorphic one. After making some choices of paths, the global symplectic topology of $E$ can be expressed in terms of an ordered collection $(L_1,\dots,L_d)$ of vanishing cycles, which are Lagrangian spheres in $M$. Denote by $\B \subset \F(M)^{opp}$ the full $A_\infty$-category, or equivalently $A_\infty$-algebra over $R$, formed by these $d$ objects (the opposite category appears here for technical reasons, having to do with our use of left modules in the body of the paper). We also have the directed subalgebra $\A \subset \B$, whose morphism spaces depend on the ordering of the objects:
\begin{equation}
e_i \A e_j = \begin{cases} e_i \B e_j & i<j, \\ \K e_i & i = j, \\ 0 & i>j. \end{cases}
\end{equation}
This has a more intrinsic meaning, being part of $\F(\pi)^{opp}$, the Fukaya category of the Lefschetz fibration as defined in \cite{seidel04}. Note that because of the Calabi-Yau (cyclic or Frobenius) nature of $\B$, the quotient $\B/\A$ is canonically isomorphic to the dual diagonal bimodule $\A^\vee[-\mathrm{dim}_{\C}(M)]$. Let $\DD$ be the curved $A_\infty$-algebra associated to the pair $(\A,\B)$.

By construction, the $L_i$ are boundaries of Lagrangian submanifolds $\Delta_i \subset E$, called Lefschetz thimbles. These are naturally objects of the wrapped Fukaya category $\W(E)$ \cite{abouzaid-seidel07} (here we are implicitly assuming that the corners of $E$ have been rounded off, so as to turn it into an exact symplectic manifold with contact type boundary). The first conjecture is:
\begin{equation} \label{eq:w-conjecture}
\parbox{36em}{
The full subcategory of $\W(E)$ with objects $(\Delta_1,\dots,\Delta_d)$ is quasi-isomorphic to the full subcategory of $\WW = modt(\DD)$ whose objects are the pullbacks of the $\A$-modules $\A e_1,\dots,\A e_d$. In particular, for each $i$ the cohomology $H(hom_{\WW}(\A e_i, \A e_i))$ is isomorphic to the wrapped Floer cohomology $HW^*(\Delta_i)$ of the corresponding Lefschetz thimble.
}
\end{equation}

\begin{example}
Take a Lefschetz fibration whose fibre $M = [-1,1] \times S^1$ is an annulus, and which has a single vanishing cycle $L_1 = \{0\} \times S^1$. In this case, $E$ is deformation equivalent to a four-dimensional ball, hence $HW^*(\Delta_1) = 0$. On the other hand, $\A = \K$, and the exact sequence \eqref{eq:b-a-sequence} obviously splits, so the category $W = H^0(\WW)$ is zero by Theorem \ref{th:main}.
\end{example}

There is a corresponding conjecture for Hochschild cohomology, which was the main subject of \cite{seidel06}:
\begin{equation}
\parbox{36em}{Up to grading-reversal, $HH_*(\DD)$ is isomorphic to the symplectic homology $SH_*(E)$. In particular, it is a symplectic invariant of $E$, independent of the Lefschetz fibration.}
\end{equation}

Take the product $D \times M$, and pick distinct cyclically ordered points $z_1,\dots,z_d \in \partial D$. Then, the disjoint union of the $K_i = \{z_i\} \times \{L_i\}$, which we denote by $K$, is a Legendrian submanifold in the boundary of $D \times M$, hence has an associated Chekanov dga, which is linear over $R$ (Chekanov's original construction \cite{chekanov99} is for Legendrian links in the standard contact $S^3$, but a generalization to the boundary of any exact symplectic manifold is envisaged as part of the general development of relative Symplectic Field Theory, see \cite{ekholm2,ekholm}). Tentatively, we'd like to propose:
\begin{equation} \label{eq:chekanov}
\parbox{36em}{The Chekanov dg algebra associated to $K$ should be quasi-isomorphic to $T(\bar{\DD}[1])^\vee$.}
\end{equation}
Here, we have used the directedness of $\A$ to write it as $\A = R \oplus \bar{\A}$, and similarly $\DD = R \oplus \bar{\DD}$. Note that Chekanov's dga is also of the form $T(\CC[1])^\vee$ for an appropriately defined $A_\infty$-algebra with curvature $\CC$. It seems natural to expect an underlying relation between $\CC$ and $\bar{\DD}$, with the caveat that the notion of equivalence for such structures is a little more tricky. 

\begin{remark}
One can think of $T(\bar{\DD}[1])^\vee$ as the endomorphism dga of the simple torsion $\DD$-module $R$, and this relates \eqref{eq:chekanov} with \eqref{eq:w-conjecture}. On the other hand, $E$ is obtained from $D \times M$ by attaching Weinstein handles to the $K_i$. This is indicative of a more general relationship between Chekanov cohomology, wrapped Floer cohomology, and Weinstein handle attachment.
\end{remark}

\subsection*{Fukaya category computations\label{subsec:fukaya}}
Take the function $(\C^*)^2 \rightarrow \C$, $(x_1,x_2) \mapsto x_1 + x_2 + (x_1x_2)^{-1}$, restrict its domain and range to suitable large compact subsets, and choose an appropriate symplectic form, so as to obtain an exact Lefschetz fibration $\pi: E \rightarrow D$. The fibre $M$ is a torus with three boundary components, and for a suitable choice of paths, the three vanishing cycles $(L_1,L_2,L_3)$ are as drawn in Figure \ref{fig:massey}.
\begin{figure}[hb]
\begin{center}
\begin{picture}(0,0)%
\includegraphics{mirrorp2.pstex}%
\end{picture}%
\setlength{\unitlength}{3947sp}%
\begingroup\makeatletter\ifx\SetFigFont\undefined%
\gdef\SetFigFont#1#2#3#4#5{%
  \reset@font\fontsize{#1}{#2pt}%
  \fontfamily{#3}\fontseries{#4}\fontshape{#5}%
  \selectfont}%
\fi\endgroup%
\begin{picture}(4436,1824)(1189,-1423)
\put(4251,-61){\makebox(0,0)[lb]{\smash{{\SetFigFont{10}{12}{\rmdefault}{\mddefault}{\itdefault}{$L_1$}%
}}}}
\put(4501,-736){\makebox(0,0)[lb]{\smash{{\SetFigFont{10}{12}{\rmdefault}{\mddefault}{\itdefault}{$L_2$}%
}}}}
\put(4801, 89){\makebox(0,0)[lb]{\smash{{\SetFigFont{10}{12}{\rmdefault}{\mddefault}{\itdefault}{$L_1'$}%
}}}}
\put(5226,-436){\makebox(0,0)[lb]{\smash{{\SetFigFont{10}{12}{\rmdefault}{\mddefault}{\itdefault}{$L_3$}%
}}}}
\end{picture}%
\caption{\label{fig:massey}}
\end{center}
\end{figure}

We will start by determining the cohomology algebra $B = H(\B)$ underlying $\B \subset \F(M)^{opp}$. It is convenient to take $\K = \C$, and to write the underlying semisimple ground ring as a group ring $R =\C[\Gamma]$, where $\Gamma = \Z/3$ (the idempotents $e_i \in R$ are then given by the three characters of the group). In these terms,
\begin{equation} \label{eq:exterior}
B = \Lambda H \semidirect \Gamma,
\end{equation}
where $\Gamma$ acts diagonally on $H = \C^3$ by cubic roots of unity, and $\Lambda H$ is the exterior algebra with the induced action. It is instructive to draw $B$ as a quiver whose vertices correspond to the idempotents $e_i \in R \subset B$, and where the arrows $i \rightarrow j$ are labeled with the spaces $e_i\bar{B}e_j$:
\begin{equation}
\xymatrix{
& \C e_3 \ar@(ul,ur)[]^-{\Lambda^3H} \ar@/^/[ddr]^-{\Lambda^2H}
 \ar@/^/[ddl]^-{H} & \\
\\
\C e_1 \ar@(l,d)[]_-{\Lambda^3H} \ar@/^/[uur]^-{\Lambda^2H} \ar@/^/[rr]^-{H} &&
\ar@/^/[ll]^-{\Lambda^2H} \C e_2 \ar@/^/[uul]^-{H} \ar@(d,r)[]_-{\Lambda^3H}
}
\end{equation}
The grading of $B$ is not the standard grading of the exterior algebra (even though the two agree mod $2$). Its precise shape depends on some choices, but one possibility is
\begin{equation}
\begin{array}{l|l}
\text{morphism space} & \
\text{degree} \\
\hline
e_1 B e_2 = H & 1 \\
e_2 B e_3 = H & 1 \\
e_3 B e_1 = H & -1 \\
e_2 B e_1 = \Lambda^2 H & 0 \\
e_3 B e_2 = \Lambda^2 H & 0 \\
e_1 B e_3 = \Lambda^2 H & 2 \\
e_k B e_k = \Lambda^0 H \oplus \Lambda^3 H & 0,1
\end{array}
\end{equation}
All this is straightforward to check, since Floer cohomology computations on $M$ can be done in a purely combinatorial fashion.

Without changing the quasi-isomorphism class of $\B$, we may assume that it is minimal ($\mu^1_\B = 0$, so it's an $A_\infty$-deformation of $B$ itself) and strictly unital. It is well-known that the classification of $A_\infty$-deformations is governed by Hochschild cohomology. In this case, only one of the relevant cohomology groups is nonzero, and
the upshot of the classification theory is as follows: there is a unique degree three homogeneous polynomial $s \in Sym^3(H^\vee)$ such that for all $h \in H$,
\begin{equation} \label{eq:vvv}
\mu^3_\B(h,h,h) = s(h).
\end{equation}
Here, the input $(h,h,h)$ is thought of as lying in $e_1 B e_2 \times e_2 B e_3 \times e_3 B e_1 = H^{\times 3}$, and the output belongs to $(e_1 B e_1)^0 = \C$. Moreover, this polynomial determines the entire $A_\infty$-structure of $\B$, up to $A_\infty$-homomorphisms whose linear term is the identity. We refer to \cite[Section 3]{seidel03b} for an exposition of the relevant algebraic deformation theory, and to \cite[Section 4]{seidel03b} for the Hochschild cohomology computation. This statement reduces the computation of $\B$ to finding $s$, hence to a finite number of unknowns. In fact, there are further constraints coming from the nontriviality of $\pi_1(M)$, which imply that $s(h) = c h_1 h_2 h_3$ for some $c \in \C$. To determine that constant, one introduces a slightly perturbed version $L_1'$ of $L_1$, and computes the Massey product by counting rectangles. This is shown in Figure \ref{fig:massey} for $h = (1,1,1)$, where precisely one relevant rectangle exists, thereby proving that $c \neq 0$ (the precise value of $c$ depends on the conventions used in the isomorphism $B \iso \Lambda H \semidirect \Gamma$, hence is irrelevant). In particular, this shows that $\B$ is not formal.

Let $\A \subset \B$ be the directed $A_\infty$-subalgebra. Because of the homogeneity of the grading, it is clear that $\A$ is a graded algebra in the ordinary sense, with vanishing higher order products. For the same reason, $\B/\A$ is an ordinary $\A$-bimodule, which means that its structure maps $\mu^{r|1|s}$ vanish for $r+s \neq 1$. In $\VV = mod(\A)$ take the simple module $R e_1$. The dual of the chain complex $hom_\VV((\B/\A)[-1] \otimes_\A Re_1, Re_1)$ is $e_1 T(\bar{\A}[1]) \otimes_R (\B/\A)[-1] \otimes_R T(\bar{\A}[1]) e_1$, which is explicitly given by
\begin{equation} \label{eq:wedge}
 H \otimes H \otimes H \xrightarrow{\begin{pmatrix} \wedge \otimes id \\ id \otimes \wedge \end{pmatrix}}
 (\Lambda^2H \otimes H) \oplus (H \otimes \Lambda^2H)
 \xrightarrow{(\wedge, -\wedge)} \Lambda^3H.
\end{equation}
The only nontrivial cohomology group of \eqref{eq:wedge} is $Sym^3H \subset H^{\otimes 3}$. Consider the canonical cocycle $\epsilon (\delta[-1] \otimes \1) \in hom_\VV((\B/\A)[-1] \otimes_\A R e_1, R e_1)$. By using the definition \eqref{eq:boundary-map} together with \eqref{eq:vvv}, one sees that its cohomology class is precisely given by $s \in Sym^3(H^\vee)$. This shows that the action of $\delta[-1]$ on the module category captures the additional information arising from the non-formality of $\B$.

\subsection*{The mirror equivalence}
$\pi: E \rightarrow D$ is the mirror of the algebraic variety $X = \CP{2}$, or more precisely (to make things basis-independent) of $\C P(H)$. A rigorous formulation of this statement is as follows. Let $V^{perf} \subset V = H^0(mod(\A))$ be the smallest thick triangulated subcategory containing the free module $\A$. Then
\begin{equation} \label{eq:mirrorp2}
\begin{parbox}{36em} {\em
$V^{perf}$ is equivalent to the derived category $D^b(X)$.
}\end{parbox}
\end{equation}
Indeed, as already pointed out in \cite{seidel00b}, $\A$ is the algebra associated to Beilinson's exceptional collection of sheaves $(\Omega^2_X(2), \Omega^1_X(1), \EuScript O_X)$. Given that, it is straightforward to construct the equivalence, which sends these sheaves to the modules ${\EuScript A} e_i$.

The bimodule $(\B/\A)[-1]$ is isomorphic to $\A^\vee[-2]$ (this is a general fact, as mentioned before, but on the other hand it's elementary to verify it in this particular case). Hence, the associated convolution function $F: V^{perf} \rightarrow V^{perf}$ is the Serre functor, up to a shift by $[-2]$. Under \eqref{eq:mirrorp2}, this corresponds to the functor of tensoring with the canonical bundle ${\EuScript K}_X$. Therefore, the natural transformation corresponding to $T$ is given by a section of ${\EuScript K}_X^{-1} \iso {\EuScript O}_X(3)$. One can use the previous computation using the simple module to show that this section is just given by the previously defined polynomial $s$ (up to a nonzero constant, which is irrelevant). The image $W^{perf}$ of $V^{perf}$ under the functor $V \rightarrow W = H^0(modt(\DD))$ is the localization along this natural transformation, which as discussed before is the derived category of coherent sheaves on the open subset $U = X \setminus Y$ for $Y = s^{-1}(0)$. For the particular $s$ which occurs here, $U \iso (\C^*)^2$. Since the object of $V$ corresponding to $\A e_1$ is a line bundle, $Hom_W(\A e_1,\A e_1)$ is just the ring of functions $\C[U]$. On the other hand, the total space $E$ of our fibration is deformation equivalent (as an exact symplectic manifold with contact type boundary) to the unit cotangent bundle of $T^2$, and this deformation also turns the Lefschetz thimble $\Delta_1$ into a cotangent fibre. It is known that in this case, the wrapped Floer cohomology is isomorphic to the homology of the based loop space, $HW^0(\Delta_1) \iso H_0(\Omega T^2;\C)$; this is in turn isomorphic to $Hom_W(\A e_1, \A e_1)$, as predicted by \eqref{eq:w-conjecture}.


\begin{thebibliography}{10}

\bibitem{abouzaid-seidel07}
M.~Abouzaid and P.~Seidel.
\newblock An open string analogue of Viterbo functoriality.
\newblock Preprint arXiv:0712.3177, 2007.

\bibitem{auroux-katzarkov-orlov04}
D.~Auroux, L.~Katzarkov, and D.~Orlov.
\newblock Mirror symmetry for weighted projective planes and their
  noncommutative deformations.
\newblock Preprint math.AG/0404281, 2004.

\bibitem{bl}
J.~Bernstein and V.~Lunts.
\newblock {\em Equivariant sheaves and functors}.
\newblock Springer Lecture Notes in Math. vol. 1578, Springer, 1994.

\bibitem{bezrukavnikov00}
R.~Bezrukavnikov.
\newblock Perverse coherent sheaves (after Deligne).
\newblock Preprint math/0005152v1, 2000.

\bibitem{chekanov99}
Yu.~Chekanov.
\newblock Differential algebras of {L}egendrian links.
\newblock {\em Invent. Math.}, 150:441--483, 2002.

\bibitem{drinfeld02}
V.~Drinfeld.
\newblock D{G} quotients of {DG} categories.
\newblock {\em J. Algebra}, 272:643--691, 2004.

\bibitem{ekholm}
T.~Ekholm.
\newblock Rational Symplectic Field Theory over $\Z_2$ for exact Lagrangian cobordisms.
\newblock Preprint math.SG/0612029.

\bibitem{ekholm2}
T.~Ekholm, J.~Etnyre, and M.~Sullivan.
\newblock The contact homology of Legendrian submanifolds in $\R^{2n+1}$.
\newblock {\em J. Differential Geom.} 71 (2005), 177--305.

\bibitem{fooo}
K.~Fukaya, Y.-G. Oh, H.~Ohta, and K.~Ono.
\newblock {\em Lagrangian intersection {F}loer theory - anomaly and obstruction.}
\newblock Book manuscript (first version, 2000; second expanded version, 2006).

\bibitem{gabriel-zisman}
P.~Gabriel, M.~Zisman.
\newblock Calculus of fractions and homotopy theory.
\newblock Springer, 1967.

\bibitem{gelfand-manin}
S.~Gelfand and Yu.~Manin.
\newblock {\em Methods of homological algebra}.
\newblock Springer, 1996.

\bibitem{keller99a}
B.~Keller.
\newblock On the cyclic homology of exact categories.
\newblock {\em J. Pure Appl. Algebra}, 136:1--56, 1999.


\bibitem{ks-notes}
M.~Kontsevich and Y.~Soibelman.
\newblock Notes on {$A_\infty$} algebras, {$A_\infty$} categories and non-commutative geometry. I.
\newblock Preprint math.RA/0606241.

\bibitem{ukraine}
V.~Lyubashenko and O.~Manzyuk.
\newblock $A_\infty$-bimodules and Serre $A_\infty$-functors.
\newblock Preprint arXiv:math/0701165, 2007.

\bibitem{orlov-branes}
D.~Orlov.
\newblock Triangulated categories of singularities and $D$-branes in Landau-Ginzburg models.
\newblock {\em Proc. Steklov Inst. Math.} 246:2004, 227--248

\bibitem{orlov-equi}
D.~Orlov.
\newblock Triangulated categories of singularities and equivalences between Landau-Ginzburg models.
\newblock Preprint math.AG/0503630, 2005.

\bibitem{polishchuk98}
A.~Polishchuk.
\newblock {M}assey and {F}ukaya products on elliptic curves.
\newblock {\em Adv. Theor. Math. Phys.}, 4:1187--1207, 2000.

\bibitem{polishchuk-zaslow98}
A.~Polishchuk and E.~Zaslow.
\newblock Categorical mirror symmetry: the elliptic curve.
\newblock {\em Adv. Theor. Math. Phys.}, 2:443--470, 1998.

\bibitem{rickard-stable}
J.~Rickard.
\newblock Derived categories and stable equivalence.
\newblock {\em J. Pure Appl. Algebra} 61:303--317, 1989.

\bibitem{seidel06}
P.~Seidel.
\newblock Symplectic homology as {H}ochschild homology.
\newblock In: {\em Proceedings of the AMS Summer Institute in Algebraic Geometry
  (Seattle, 2005)}, Amer. Math. Soc., in press.

\bibitem{seidel00b}
P.~Seidel.
\newblock More about vanishing cycles and mutation.
\newblock In K.~Fukaya, Y.-G. Oh, K.~Ono, and G.~Tian, editors, {\em
  {S}ymplectic {G}eometry and {M}irror {S}ymmetry ({P}roceedings of the 4th
  {KIAS} Annual International Conference)}, pages 429--465. World Scientific,
  2001.

\bibitem{seidel03b}
P.~Seidel.
\newblock Homological mirror symmetry for the quartic surface.
\newblock Preprint math.SG/0310414, 2003.

\bibitem{seidel04}
P.~Seidel.
\newblock {\em {F}ukaya categories and {P}icard-{L}efschetz theory}.
\newblock European Math. Soc., to appear.

\bibitem{tradler01}
T.~Tradler.
\newblock Infinity-inner-products on {A}-infinity-algebras.
\newblock Preprint math.\-AT/0108027.

\bibitem{weibel}
C.~Weibel.
\newblock {\em An introduction to homological algebra}.
\newblock Cambridge Univ. Press, 1994.
\end{thebibliography}
\end{document}